%% file: main_arxiv.tex
\documentclass[review]{elsarticle}

\journal{~}

\usepackage{etoolbox}
\makeatletter
\patchcmd{\ps@pprintTitle}% <cmd>
{Preprint submitted to}% <search>
{~}% <replace>
{}{}% <succes><failure>
\makeatother

\usepackage[top=2.5cm,bottom=5cm,left=3cm,right=3cm]{geometry}

\usepackage{graphics}
\usepackage{bm}
\usepackage{amsmath}
\usepackage{amsfonts}
\usepackage{xcolor}
\usepackage{manfnt}
\usepackage[utf8]{inputenc}
\usepackage[T1]{fontenc}
\usepackage[ngerman,english]{babel}
\usepackage{booktabs}

\newcommand{\DD}{Data-Driven}
\newcommand{\Dd}{Data-driven}
\newcommand{\dd}{data-driven}

\newcommand{\DDCM}{\DD\ Computational Mechanics}

\input{commands_r1}

\begin{document}

\input{0front}

%\linenumbers

\tableofcontents

\newpage

\input{1introduction}
\input{2mechanics}

\input{3programs}
\input{4structure}
\input{5examples}
\input{6conclusion}
\input{7back}

%\section*{References}
\bibliographystyle{elsarticle-num}
\bibliography{bib}

\end{document}

%% file: commands_r1.tex
\newcommand{\field}{\mathbb}

\newcommand{\reals}{\field{R}}

\newcommand{\euclideans}{\field{E}}

\newcommand{\spaceX}{\mathcal{X}}

\newcommand{\spaceW}{\mathcal{W}}

\newcommand{\spaceWdual}{\mathcal{W}^{*}}
\newcommand{\xthetat}{\bm{x}(\bm{\theta};t)}
\newcommand{\xthetaz}{\bm{x}(\bm{\theta};0)}
\newcommand{\dlangle}{\langle\!\langle}
\newcommand{\drangle}{\rangle\!\rangle}

\newcommand{\cube}{\mbox{\mancube}}

\newcommand{\dvolGd}{\sqrt{\det[\bmetric_d(\bm{\theta};0)]}\,\mathrm{d}\,\cube}

\newcommand{\bmetric}{\bm{G}}

\newcommand{\ie}{\textit{i.e.,}~}
\newcommand{\eg}{\textit{e.g.,}~}
\newcommand{\ia}{\textit{inter alia}}
\newcommand{\etal}{\textit{et al.}~}

\newcommand{\etc}{\textit{etc}}

\newcommand{\uii}{\bm{i}^1}
\newcommand{\ujj}{\bm{i}^2}
\newcommand{\ukk}{\bm{i}^3}

\newcommand{\zero}{\bm{0}}
\newcommand{\coor}{\bm{q}}
\newcommand{\weight}{\bm{C}}
\newcommand{\strain}{\bm{e}}
\newcommand{\stress}{\bm{s}}

\newcommand{\bxi}{\bm{\xi}}
\newcommand{\force}{\bm{f}}
\newcommand{\lagrangean}{\mathcal{L}}
\newcommand{\cost}{\mathcal{J}}
\newcommand{\cmanifold}{\mathcal{Z}}
\newcommand{\datacmanifold}{\widetilde{\cmanifold}}
\newcommand{\reccmanifold}{\check{\cmanifold}}
\newcommand{\blambda}{\bm{\lambda}}
\newcommand{\bmu}{\bm{\mu}}
\newcommand{\bnu}{\bm{\nu}}
\newcommand{\bchi}{\bm{\chi}}
\newcommand{\bphi}{\bm{\varphi}}
\newcommand{\btheta}{\bm{\theta}}
\newcommand{\pd}[2]{\partial_{#1}{#2}}

\newcommand{\ba}{\bm{a}}
\newcommand{\bb}{\bm{b}}
\newcommand{\bg}{\bm{g}}
\newcommand{\bx}{\bm{x}}

\newcommand{\bn}{\bm{n}}
\newcommand{\gl}{\bm{E}_{\flat}}
\newcommand{\rotation}{\bm{\Lambda}}
\newcommand{\KKT}{\bm{S}}
\newcommand{\refe}{\mathrm{ref}}

\newcommand{\ofcoor}{(\coor)}
\newcommand{\bgamma}{\bm{\gamma}}
\newcommand{\bomega}{\bm{\omega}}
\newcommand{\coorz}{\bm{\varphi}_0}
\newcommand{\coori}{\bm{d}_1}
\newcommand{\coorj}{\bm{d}_2}
\newcommand{\coork}{\bm{d}_3}
\newcommand{\bforce}{\bm{n}}
\newcommand{\bmoment}{\bm{m}}

\newcommand{\bE}{\bm{E}}

\newcommand{\bI}{\bm{I}}
\newcommand{\bB}{\bm{B}}
\newcommand{\bh}{\bm{h}}
\newcommand{\bH}{\bm{H}}
\newcommand{\bU}{\bm{U}}
\newcommand{\bUbar}{\bar{\bU}}
\newcommand{\bUbre}{\breve{\bU}}
\newcommand{\bV}{\bm{V}}
\newcommand{\bW}{\bm{W}}

\newcommand{\coorihat}{\widehat{\coori}}
\newcommand{\coorjhat}{\widehat{\coorj}}
\newcommand{\coorkhat}{\widehat{\coork}}

\newcommand{\bbtwo}{\bb_2}
\newcommand{\bbtwohat}{\widehat{\bbtwo}}

\newcommand{\bai}{\ba_1}
\newcommand{\baj}{\ba_2}
\newcommand{\bak}{\ba_3}
\newcommand{\baihat}{\widehat{\bai}}
\newcommand{\bajhat}{\widehat{\baj}}
\newcommand{\bakhat}{\widehat{\bak}}

\newcommand\argvar
  {(\coor, \strain, \stress, \blambda, \bmu, \bnu; \straindata, \stressdata)}
\newcommand{\argvarrec}{(\strainrec, \stressrec, \coor, \strain, \stress, \blambda, \bmu, \bnu, \bxi )}
\newcommand{\cons}{\bm{h}(\coor)}
\newcommand{\nullproj}{\bm{N}(\coor)}
\newcommand{\nullprojT}{\bm{N}(\coor)^T}

\newcommand{\straincoor}{\strain(\coor)}

\newcommand{\straindata}{\widetilde{\strain}}
\newcommand{\strainrec}{\check{\strain}}
\newcommand{\iweight}{\weight^{-1}}
\newcommand{\stressdata}{\widetilde{\stress}}
\newcommand{\stressrec}{\check{\stress}}

\newcommand\KKTcoorcoor{\KKT_{\coor\coor}(\coor, \stress, \blambda, \bmu, \bnu)}
\newcommand\KKTstresscoor{\KKT_{\stress\coor}(\coor, \bmu)}
\newcommand\KKTmucoor{\KKT_{\bmu\coor}(\coor, \stress)}
\newcommand\KKTmustress{\KKT_{\bmu\stress}(\coor)}

\newcommand{\bHcoor}{\bH\ofcoor}
\newcommand{\bBcoor}{\bB\ofcoor}

\newcommand{\bHcoorT}{\bH\ofcoor^T}
\newcommand{\bBcoorT}{\bB\ofcoor^T}

\newcommand\nstrain{n_{\strain}}

\newcommand\fixNLP{NLP$(\straindata, \stressdata)$}
\newcommand\MCS[1][?]{\textcolor{violet}{[MCS: #1]}}
\newcommand\sdot{^T}
\newcommand\sfix{_{\mathrm{fix}}}
\newcommand\sone{}

\newcommand\rdot{{}\cdot{}}
\newcommand\tfig[1][48]{\includegraphics[width=0.#1\textwidth]}
\newcommand\dfig[2][48]{\tfig[#1]{noise00-fig_#2}}
\newcommand\sfig[2][48]{\tfig[#1]{noise10-fig_#2}}

\newcommand{\arclen}{\sigma}
\newcommand{\pdarclen}{\partial_\arclen}

\newcommand{\Boper}{\mathfrak{B}}
\newcommand{\Uoper}{\mathfrak{U}}

\newcommand{\Doper}{\mathfrak{D}}

%% file: 0front.tex
\begin{frontmatter}

\title{\ A Framework for \DDCM\\ Based on Nonlinear Optimization}

\author[add1]{Cristian Guillermo Gebhardt\corref{cor1}}
\ead{c.gebhardt@isd.uni-hannover.de}
\cortext[cor1]{Corresponding author}%, Tel.: +49-511-762-2859, Fax: +49-511-762-2236.}
\author[add2]{Dominik Schillinger}
\ead{schillinger@ibnm.uni-hannover.de}
\author[add3]{Marc Christian Steinbach}
\ead{mcs@ifam.uni-hannover.de}
\author[add1]{Raimund Rolfes}
\ead{r.rolfes@isd.uni-hannover.de}
\address[add1]{Leibniz Universität Hannover, Institute of Structural Analysis - ForWind Hannover,\\ Appelstraße 9A, 30167 Hannover, Germany}
\address[add2]{Leibniz Universität Hannover, Institute of Mechanics and Computational Mechanics,\\ Appelstraße 9A, 30167 Hannover, Germany}
\address[add3]{Leibniz Universität Hannover, Institute of Applied Mathematics,\\ Welfengarten 1, 30167 Hannover, Germany}

\begin{abstract}
\DDCM~is a novel computing paradigm that enables the transition from standard data-starved approaches to modern data-rich approaches. At this early stage of development, one can distinguish two mainstream directions. The first one relies on a discrete-continuous optimization problem and seeks to assign to each material point a point in the phase space that satisfies compatibility and equilibrium, while being closest to the data set provided. The second one is a \dd~inverse approach that seeks to reconstruct a constitutive manifold from data sets by manifold learning techniques, relying on a well-defined functional structure of the underlying constitutive law. In this work, we propose a third route that combines the strengths of the two existing directions and mitigates some of their weaknesses. This is achieved by the formulation of an approximate nonlinear optimization problem, which can be robustly solved, is computationally efficient, and does not rely on any special functional structure of the reconstructed constitutive manifold. Additional benefits include the natural incorporation of kinematic constraints and the possibility to operate with implicitly defined stress-strain relations. We discuss important mathematical aspects of our approach for a \dd~truss element and investigate its key numerical behavior for a \dd~beam element that makes use of all components of our methodology.
\end{abstract}


\end{frontmatter}

%% file: 1introduction.tex
\section{Introduction}

In classical solid mechanics, a set of constraints and conservation laws, such as the equations of  compatibility and equilibrium, describe the deformation of material bodies \cite{Oden2011}. These equations result from fundamental physical principles and hold for any material body irrespective of its material properties. To close the associated boundary value problem, constraints and conservation laws must be accompanied by constitutive equations that relate the quantities of the phase space, such as stresses and strains. In classical solid mechanics, constitutive equations are provided by phenomenological material models that can be calibrated via observational data. The design, mathematical formulation and numerical solution of material models has been an important research field since the beginnings of solid mechanics \cite{Simo2006,De2011} and remains the subject of extensive ongoing research to date.

In contrast to balance laws that have axiomatic character, material models are empirical and therefore constitute a source of error and uncertainty, especially when the material behavior is complex. Given that experimental data on the constitutive behavior is scarce, material modeling seems to be an inevitable step. Over the last few years, however, materials science has been undergoing a remarkable transition from a data-starved to a data-rich field, with an increasing number of scenarios where an abundance of data is available to characterize the constitutive behavior of materials. This is mainly due to technological advances in the field of experimental measurements, data storage and data processing among others.
In this context, \DDCM, a new computing paradigm originally initiated by Kirchdoerfer and Ortiz \cite{Kirchdoerfer2016}, is currently emerging. Its idea is to reformulate classical boundary value problems of elasticity and inelasticity in such a way that empirical material models are replaced by experimental material data described in the phase space. Its overarching goal is to eliminate the modeling error and uncertainty of phenomenological material models and instead directly exploit the available wealth of experimental data in its entirety. However, a new source of error appears that is related to the measurements in themselves and to the measuring chain. Up to this point, it is not clear, which of the error sources is the most harmful.   
\ifcase1% Q
\MCS[Do Kirchdoerfer, Ortiz etc.\ claim that direct use of data reduces modeling
errors and uncertainty? If yes, why? I would argue that the opposite is true!]
A CGG: Yes, they do, but we don't ;)
\fi

At this early stage of development, one can differentiate two basic approaches in \DDCM. Following Kirchdoerfer and Ortiz \cite{Kirchdoerfer2016}, one class of methods results in
a \dd~solver that seeks to assign to each material point a point in the phase space, which besides satisfying compatibility and equilibrium is closest to the data set provided. Its formulation is based on a discrete-continuous optimization problem that minimizes the distance between the material data set and the subspace of compatible strain fields and stress fields in equilibrium. First attempts exist that extend this approach beyond linear elasticity, for instance to geometrically nonlinear elasticity \cite{Nguyen2018}, general elasticity \cite{Conti2018}, elastodynamics \cite{Kirchdoerfer2018}, and inelasticity \cite{Eggersmann2019}. The current discrete-continuous approach, however, is computationally expensive and missing robustness in certain situations.
In particular, it exhibits strong sensitivity to scattering of the data set, and solution by meta-heuristic optimization techniques leads to relatively poor convergence. First attempts have been published that try to mitigate these drawbacks, \eg a maximum entropy scheme that increases robustness with respect to outliers~\cite{Kirchdoerfer2017}, or the mathematically well-behaved formulation of the discrete-continuous optimization problem as a computationally tractable mixed-integer quadratic optimization problem \cite{Kanno2019}.

The second class of methods follows the \dd~inverse approach of Cueto, Chinesta and collaborators \cite{Ibanez2018a} who seek to reconstruct a constitutive manifold from data, using manifold learning methods. In the case of elasticity, the goal is to use data to identify a suitable approximation of the strain energy density functional, whose first and second derivatives result in the stress tensor and the elastic tangent, respectively. In a broader context \cite{Ibanez2017}, it is proposed to identify the locally linear behavior and for non-convergent cases, it is proposed to find the intersection between the equilibrium and constitutive manifolds. And recently, this approach has been successfully applied to the setting of ``General Equation for Non-Equilibrium Reversible-Irreversible Coupling'' \cite{Gonzalez2019}. Although this approach has a number of advantages, such as straightforward reassurance of thermodynamic consistency, the transfer of data-intensive computations in an off-line step, the potential for nonintrusive implementation in standard codes, it also entails a number of significant limitations. Most importantly, the method relies in the assumption of constitutive manifolds with a special functional structure and is thus limited to the explicit definition of stress, at least in its incremental form.
% Dominik: we need more citations in this paragraph. There exists significant work in this direction in France, that unfortunately starts to deviate significantly from the basic approaches. We should think about how we take this into account.

The goal of the current work is to explore a synergistic compromise between these two classes of methods that combines their strengths and mitigates some of their main weaknesses. Our main focus is the formulation of an approximate nonlinear optimization problem. On the one hand, it improves computational efficiency and robustness with respect to the current discrete-continuous optimization approach of the first class of \dd~solvers. In particular, our approximate nonlinear optimization problem can be solved locally with standard Newton type nonlinear optimization methods, without the need to resort to more exotic options such as meta-heuristic methods that lack well-matured mathematical foundations when compared to gradient-based methods. On the other hand, it does not rely on special assumptions of the functional structure of the reconstructed constitutive manifold as the second class of data-riven solvers that potentially re-introduce modeling errors and uncertainty. Additional benefits of our optimization approach include the natural incorporation of kinematic constraints and the possibility to operate with implicitly defined stress-strain relations, which enlarges the range of material behavior that can be addressed. As our primary goal is a proof of concept for our new approach, we will use artificial rather than real measurements when required. We showcase the advantages of our approach for the case of a \dd~geometrically exact beam element that makes use of all components of our methodology.
% Dominik: please feel free to correct the details.

The structure of the article is as follows: Section~2 provides a concise review of some pertinent elements of computational nonlinear solid mechanics and fixes our notation and terminology that we will use in the remainder of the article. In Section~3, we first describe the full discrete-continuous optimization problem in its global format. We then describe in detail our new approach based on an approximate nonlinear optimization problem. In particular, we discuss the approximated implicitly defined constitutive manifold, the associated Lagrangian functional and the Karush–Kuhn–Tucker (KKT) conditions, the linearization of such KKT conditions and the resulting KKT matrix in explicit format. We note that the formulations presented in Section~3 are general and only assume a discretization of the problem in the context of the finite element method. In Section~4, we first illustrate our optimization approach for a simple truss element and then proceed to the more complex case of the geometrically exact beam element. Section~5 presents numerical examples for the geometrically exact beam that demonstrate the advantages of our approach. We close with a summary and the main conclusions in Section~6.

%% file: 2mechanics.tex
\section{Framework of computational nonlinear solid mechanics}

% Mechanical framework

\subsection{Continuous setting}

\begin{figure}[tp]
	\centering
	\includegraphics[trim=30 60 30 80,clip,width=0.7\textwidth]{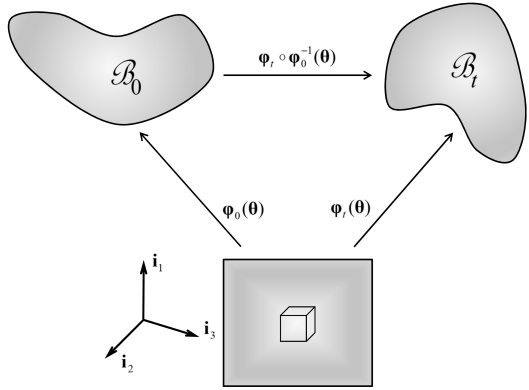}
	\caption{The continuum body: evolution among configurations through the regular motion $\bphi_t\circ\bphi_0^{-1}$.}
	\label{fig:continuum}
\end{figure}

\noindent Let us consider a continuum body characterized by two reference
sets denoted by $\mathcal{B}_{0}$, the original configuration, and
$\mathcal{B}_{t}=\bm{\varphi}_{t} \circ \bm{\varphi}^{-1}_{0} (\mathcal{B}_{0})$, the current configuration, both open sets of
$\euclideans^{3}$ ($\reals^3$ with the standard Euclidean structure), see Figure  \ref{fig:continuum}. % at time instants $t_{\xi} = 0$ and $t_{\xi} = t$, respectively.
Here $\bm{\varphi}_t \circ \bm{\varphi}^{-1}_{0}$ is a smooth regular motion from the original configuration to the current one, \ie its inverse and derivatives are well defined everywhere, and the symbol $\circ$ denotes a suitable composition rule. The parameter $t$ is a smooth parameter intended for the indexation of all possible successive configurations. Here, a chart is given by the pair $(\bm{\theta},\bm{\varphi}_t)$, with $\bm{\theta}$ a subset of $\reals^3$ and the mapping function $\bm{\varphi}_t\colon\bm{\theta} \mapsto \xthetat$.
The configuration of the system is described by the vector field
$\xthetat\in\spaceX\subseteq\euclideans^{3}$, and the static problem can be thoroughly formulated by means of the following (non-variational) set of equations:
%behavior within the bounded time interval $[0,T]\subset\reals_{\geq 0}$ can be weakly formulated as
%
\begin{subequations}
\begin{align}
\label{eq1a} \zero &= \bm{E}_{\flat}-\bm{E}_{\flat}(\bm{x}),\\
\label{eq1b} \zero &= \bm{f}^{\textrm{int}}(\bm{x},\bm{S}^{\sharp})+\bm{H}(\bm{x})^T\bm{\chi}-\bm{f}^{\textrm{ext}},\\
\label{eq1c} \zero &= \bm{h}(\bm{x}).
\end{align}
\end{subequations}
Here $\flat$ and $\sharp$ indicate that a rank-$n$ tensor is $n$ times covariant and $n$ times contravariant, respectively.
% Dominik: This definition needs to be known to understand the following equations, mentioning it a page later might be too late...

In the first equation \eqref{eq1a}, the displacement-based Green-Lagrange strain tensor, an element of $\mathcal{E}:=\{\bm{E}_{\flat}(\bm{x}) \in T^{*}_{\xthetaz}\spaceX \times T^{*}_{\xthetaz}\spaceX \mid \textrm{skew}(\bm{E}_{\flat}(\bm{x}))=\bm{0}\}$ with $\textrm{skew}(\cdot)$ indicating the skew-symmetric part of the tensor considered, is given in the curvilinear setting by
%
% Dominik: I removed superscript "disp" for eps. Are later version of eps that do not carry that superscript ever not displacement-based?
\begin{equation}
\bm{E}_{\flat}(\bm{\theta};t) = \frac{1}{2}\left(\bmetric(\bm{\theta}; t)-\bmetric(\bm{\theta}; 0)\right).
\end{equation}
The pullback of the metric tensor at $\mathcal{B}_t$ through the regular motion $\bm{\varphi}_t \circ \bm{\varphi}^{-1}_{0}$, \ie~$\bmetric(\bm{\theta}; t)\colon T^{*}_{\xthetaz}\spaceX \times T^{*}_{\xthetaz}\spaceX \to \reals_{\geq 0}$, is defined as
%
%\newpage
%
\begin{equation}
\bmetric(\bm{\theta}; t) = \delta_{ij}\pd{\theta^a}{x^i(\bm{\theta}; t)}\pd{\theta^b}{x^j(\bm{\theta}; t)}\, \bm{g}^a(\bm{\theta}; 0) \otimes \bm{g}^b(\bm{\theta}; 0),
\label{eq:metric}
\end{equation}
where $\delta_{ij}$ are the components of the Euclidean metric tensor and $\otimes$ denotes the outer product. The elements $\bm{g}_{a}\in T_{\xthetaz}\spaceX$ of the contravariant basis are defined as $\bm{g}_{a}=\pd{\theta^{a}}{x^{i}}\bm{i}_{i}$, where $\bm{i}_{i}$ with $i$ from $1$ to~$3$ contains the elements of the standard orthonormal basis in $\euclideans^{3}$, \ie the space of column vectors. The elements $\bm{g}^{a}\in T_{\xthetaz}^{*}\spaceX$ of the covariant basis are defined such that $\langle\bm{g}^{b},\bm{g}_{a}\rangle =\delta_{b}^{a}$. % where $\bm{g}_{a}$ are the elements of the contravariant basis.
The admissible variation of the displacement-based Green-Lagrange strain tensor is denoted by $\delta \bm{E}_{\flat}(\bm{x})$ and, as stated previously, this is intrinsically related to $\delta\bm{x}$.

In the second equation \eqref{eq1b}, the so-called balance equation, we use the property
\begin{equation}
\int_{\mathcal{B}_{0}}
\langle \delta \bm{x}, \bm{f}^{\mathrm{int}}(\bm{x},\bm{S}^{\sharp})\rangle
\, \mathrm{d}\mathcal{B}_{0}
 =
\int_{\mathcal{B}_{0}}
\dlangle\delta\bm{E}_{\flat}(\bm{x}),\bm{S}^{\sharp}\drangle
\, \mathrm{d}\mathcal{B}_{0}
,
\end{equation}
to define $\bm{f}^{\mathrm{int}} \in T_{\xthetat}^{*}\spaceX$, the vector density of internal forces. The angle brackets $\langle\rdot,\rdot\rangle\colon\spaceW\times\spaceWdual \to \reals$ denote an appropriate dual pairing, in which $\spaceW$ is a vector space (whose elements are called vectors) and $\spaceWdual$ is its dual space (whose elements are called covectors or one-forms), and the double brackets $\dlangle\rdot,\rdot\drangle\colon (\spaceW\times\spaceW)\times(\spaceWdual\times\spaceWdual) \to \reals$ represent an appropriate double dual pairing. $\delta\bm{x}\in T_{\xthetat}\spaceX$ denotes any admissible variation of the configuration vector.
Moreover, $\bm{H}(\bm{x})\in L(T_{\xthetat}\spaceX, \reals^m)$ is the Jacobian matrix of $\bm{h}(\bm{x})=\bm{0}\in\reals^m$, a finite-dimensional field of integrable constraints specified in the third equation \eqref{eq1c}. These constraints can be used for instance to impose boundary conditions. $\bm{\chi}\in\reals^m$ represents the corresponding Lagrange multipliers, and $\bm{f}^{\textrm{ext}}\in T_{\xthetat}^{*}\spaceX$ denotes the vector density of external forces.

Furthermore, it is worthwhile to note that the Green-Lagrange strain tensor, $\bm{E}_{\flat}$, %:=\{\bm{E}_{\flat} \in T^{*}_{\xthetaz}\spaceX \times T^{*}_{\xthetaz}\spaceX \mid \textrm{skew}(\bm{E}_{\flat})=\bm{0}\}$,
% Dominik: Already defined above...
is to be separately accounted for as a primal field.
 The specific strain measure chosen is accompanied by the second Piola-Kirchhoff stress tensor, an element of $\mathcal{S}:=\{\bm{S}^{\sharp} \in T_{\xthetaz}\spaceX \times T_{\xthetaz}\spaceX \mid \textrm{skew}(\bm{S}^{\sharp}) = \bm{0}\}$, which in the case of general hyperelastic materials is given by $\bm{S}^{\sharp}=\pd{\bm{E}_{\flat}}{\Psi^{\mathrm{int}}(\bm{E}_{\flat})}$, where $\Psi^{\mathrm{int}}(\bm{E}_{\flat})$ denotes the associated strain energy. We recall that a three-field variational principle, for instance, requires this particular functional structure. In the context of the present work, we do not rely on the standard functional structure between the strain and stress tensors. If one instead considers the existence of a strain energy function and strongly imposes the compatibility condition, which can be done by direct substitution of the strain tensor (\ie replacing $\bm{E}_\flat$ by $\bm{E}_{\flat}=\bm{E}_{\flat}(\bm{x})$ everywhere when needed), then, the governing equations can be reverted to a fully variational set of equations.

\subsection{Discrete setting}

\noindent By using, for instance, isoparametric finite elements, the governing equations for a single element can be discretely approximated as
\begin{subequations}
	\begin{align}
	\zero &= \int_{\cube}\bm{\Phi}^T_{\stress}\left(\strain^h-\strain(\bx^h)\right)\dvolGd,\\
	\zero &= \int_{\cube}\bm{\Phi}^T_{\bx}\left(\Boper(\bx^h)^T\stress^h+\bH(\bx^h)^T\bchi^h-\force^{\mathrm{ext},h}\right)\dvolGd,\\
	\zero &= \int_{\cube}\bm{\Phi}^T_{\bchi}\,\bm{h}(\bx^h)\dvolGd.
	\end{align}
\end{subequations}
The superscript ``$h$'' denotes that the field is discretized in space by means of the interpolation of nodal or elemental discrete variables and the matrices $\bm{\Phi}_{\bx}$, $\bm{\Phi}_{\stress}$ and $\bm{\Phi}_{\bchi}$ contain the admissible test functions chosen. In addition, the strain and stress tensor fields are represented by an adequate notation that takes advantage of their symmetries, \ie $\bm{E}_\flat\leadsto\strain$ and $\bm{S}^\sharp\leadsto\stress$. As a direct consequence, the discrete version of $\force^{\mathrm{int}}(\bx,\bm{S}^\sharp)$ can be written as $\Boper(\bx^h)^T\stress^h$ and the discrete version of $\bm{E}_{\flat}-\bm{E}_{\flat}(\bm{x})$ as $\strain^h-\strain(\bx^h)$. Moreover, $\Boper(\bx^h)\leadsto\bB(\coor)$, where the approximated displacement field is computed in terms of a finite dimensional set of generalized coordinates, \ie $\bx^h(\btheta,\coor)=\bm{\Phi}_{\bx}(\btheta)\coor$. $\btheta$ spans the %one-, two- or three-dimensional
master domain symbolized as $\cube$, and $\dvolGd$ is the discrete version of the volume element. Finally, the integrals are computed by means of a given integration scheme, for instance Gauss quadrature, yielding the final format of the discrete governing equations that will be used in the remainder of this article.

%% file: 3programs.tex
\section{Optimization problems for \DDCM}

In the context of \DDCM, the simplest scalar cost function
to be minimized can be defined as \cite{Kirchdoerfer2016}
\begin{equation}
  \label{eq:cost}
  \cost(\straindata, \stressdata, \strain, \stress)
  =
  \frac12 \|\strain - \straindata\|^2_{\weight} +
  \frac12 \|\stress - \stressdata\|^2_{\iweight}
  \, .
\end{equation}
Here $(\strain, \stress) \in \cmanifold$ comprises continuous
strain and stress variables $\strain$ and $\stress$, respectively,
the given data set $\datacmanifold$ contains
finitely many strain and stress measurements
$(\straindata, \stressdata) \in \datacmanifold$,
$\weight \in \reals^{\nstrain \times \nstrain}$
(with $\nstrain = \dim(\strain) = \dim(\stress)$)
is a symmetric positive-definite weight matrix
with inverse $\iweight$,
and $\|\rdot\|_{\weight}$ and $\|\rdot\|_{\iweight}$ are norms derived from the inner product.
At this point, we do not specify the dimension $\nstrain$
because it depends strongly on the model considered.
For instance, we have $\nstrain = 6$ for a Reissner-Sim\'{o} beam
\cite{Romero2004} and $\nstrain = 4$ for a Kirchhoff rod \cite{Romero2019},
and therefore, we leave this detail out of the current discussion.
The cost function \eqref{eq:cost} has to be minimized
under the following constraints\footnote{%
  As this section focuses on the overall optimization problem,
  we skip details on domain integration,
  quadrature scheme chosen, \etc.
  Further details on such aspects can be consulted elsewhere.}:
\textit{i}) the compatibility equation that enforces the equivalence
between strain variables and displacement-based strains,
\begin{equation}
  \strain - \straincoor = \zero,
\end{equation}
in which $\coor \in Q \subset \reals^{m + n}$
is the vector of generalized coordinates
employed to describe the system kinematics,
with $Q$ standing for the configuration manifold;
\textit{ii}) the balance equation that establishes the static equilibrium,
\begin{equation}
  \bBcoorT \stress + \bHcoorT \bchi - \force = \zero,
\end{equation}
in which
$\bBcoor = \pd\coor \strain(\coor) \in L(T_{\coor} Q, \reals^{\nstrain})$
is the Jacobian matrix of the displacement-based strains,
$\bHcoor = \pd\coor \bh(\coor) \in L(T_{\coor} Q, \reals^m)$
is the Jacobian matrix of the kinematic constraints,
$\bchi \in \reals^m$
is the corresponding vector of Lagrange multipliers,
and $\force\in T_{\coor}^*Q$ represents the vector of generalized external loads
(the superscript ``$\mathrm{ext}$'' has been dropped);
and, \textit{iii}) the kinematic constraints
\begin{equation}
  \cons = \zero,
\end{equation}
a finite set of integrable restrictions that belongs to $\reals^m$.
In the finite element context,
the operator $\bBcoor$ and the Jacobi matrix $\bHcoor$ are typically linear in $\coor$,
which simplifies the computation of related higher-order derivatives enormously.
We will see this later on in Section~4.
A counterexample to this is, for instance,
a geometrically exact beam element in which the rotations are parametrized
with the Cartesian rotation vector \cite{Cardona1988}.
To avoid problems caused by overdetermination
and singular KKT matrices in the subsequent optimization problems,
we eliminate the Lagrange multipliers from the balance equation.
This goal can be achieved by several techniques,
such as projection onto a space orthogonal to the constraint forces,
reduction by means of coordinate partitioning, or redefinition
of the derivative on the configuration manifold (connection),
see for instance \cite{Heard2006}.
Here, we choose the first approach,
known as the null-space approach,
which requires a null-space basis matrix
$\nullproj \in L(\reals^n, T_{\coor} Q)$ for
$\ker(\bHcoor) = \{ \bn \in T_{\coor} Q \mid \bHcoor \bn = \zero \in \reals^m \}$
with $n = \dim(Q) - m$, the system's number of degrees of freedom,
and $\mathrm{rank}(\nullproj) = \dim(\ker(\bHcoor)) = n$,
such that
\begin{equation}
  \bHcoor \nullproj = \zero.
\end{equation}
Then the balance equation adopts the form
\begin{equation}
  \nullprojT (\bBcoorT \stress - \force) = \zero
\end{equation}
and further on, the (unique) vector of Lagrange multipliers is
\begin{equation}
  \bchi = -(\bHcoor \bHcoorT)^{-1} \bHcoor (\bBcoorT \stress - \force),
\end{equation}
and the unique orthogonal null-space projector can be expressed as
\begin{equation}
  \bm P(\coor) = \bI - \bHcoorT(\bHcoor \bHcoorT)^{-1} \bHcoor.
\end{equation}
Additionally, an infinite number of non-orthogonal
(oblique) null-space projectors exist.
In practice, it is often more efficient to obtain the vectors
$\nullprojT (\bBcoorT \stress - \force)$ and $\bchi$
from an implicit representation of $\nullproj$
than forming $\bm P(\coor)$ explicitly.
Details can be found in standard textbooks
on numerical linear algebra or on optimization,
such as \cite{Golub-van_Loan:1989,Nocedal-Wright:2006}.
In the context of slender structures, the reader may refer to
\cite{Gebhardt2017,Hente2019} and references therein
for the efficient calculation of the null-space basis.

\subsection{A discrete-continuous nonlinear optimization problem}

Employing directly the strain and stress measurements,
a discrete-continuous nonlinear optimization problem
(brief\/ly called a DCNLP here) can be stated as
\begin{equation}
  \label{eq:dcnlp}
  \begin{split}
    \min_{(\straindata,\stressdata,\coor,\strain,\stress)} \quad &
    \frac12 \|\strain - \straindata\|^2_{\weight} +
    \frac12 \|\stress - \stressdata\|^2_{\iweight} \\
    \text{subject to} \quad
    &\strain - \straincoor =\zero, \\
    &\nullprojT (\bBcoorT \stress - \force) = \zero, \\
    &\cons = \zero.
  \end{split}
\end{equation}
Note that the discrete variables $(\straindata,\stressdata)$
appear only in the cost function.
In \eqref{eq:dcnlp} there is no mathematical structure
that relates these discrete variables to each other;
a mathematically well-behaved formulation with a binary selection variable
for each pair $(\straindata,\stressdata) \in \datacmanifold$
would be called a mixed-integer nonlinear optimization problem (MINLP).
For fixed $(\straindata,\stressdata) \in \datacmanifold$,
the problem becomes a smooth nonlinear optimization problem (NLP),
referred to as \fixNLP, with associated Lagrangian function
\begin{equation}
  \begin{split}
    \lagrangean\sfix\argvar
    &=
    \frac12 \|\strain - \straindata\|^2_{\weight} +
    \frac12 \|\stress - \stressdata\|^2_{\iweight} \\
    &+ \blambda \sdot (\strain - \straincoor) \\
    &+ (\nullproj \bmu) \sdot (\bBcoorT \stress - \force) \\
    &+ \bnu \sdot \cons,
  \end{split}
\end{equation}
where $\blambda \in \reals^{\nstrain}$ are Lagrange multipliers
of the compatibility equation,
$\bmu \in \reals^n$ are Lagrange multipliers
of the balance equation premultiplied by the null-space basis matrix, and
$\bnu \in \reals^m$ are Lagrange multipliers of kinematic constraints.
For a fixed measurement pair $(\straindata,\stressdata)$,
any solution of \fixNLP\ provides
a set of values $(\coor,\strain,\stress)$
that locally minimizes the cost function under the given constraints.
To derive the corresponding first-order optimality conditions,
we calculate the variation of $\lagrangean\sfix$ as
\begin{equation}
  \begin{split}
    \delta\lagrangean\sfix\argvar
    &=
    \delta\coor \sdot (-\bBcoorT \blambda +
    (\pd{\coor} [\nullprojT (\bBcoorT \stress - \force)])^T \bmu +
    \bHcoorT \bnu) \\
    &+ \delta\strain \sdot (\weight (\strain - \straindata) + \blambda) \\
    &+ \delta\stress \sdot
        (\iweight (\stress - \stressdata) + \bBcoor \nullproj \bmu) \\
    &+ \delta\blambda \sdot (\strain - \straincoor) \\
    &+ \delta\bmu \sdot \nullprojT (\bBcoorT \stress - \force) \\
    &+ \delta\bnu \sdot \cons.
  \end{split}
\end{equation}
This variation has to vanish for any choice of the varied quantities,
and hence we obtain the following primal-dual system of equations,
which are known as the KKT conditions in optimization:
\begin{subequations}
  \begin{align}
    &\delta\coor:&
    -\bBcoorT \blambda +
    (\pd{\coor} [\nullprojT (\bBcoorT \stress - \force)])^T \bmu +
    \bHcoorT \bnu &= \zero, \\
    &\delta\strain:&
    \weight (\strain - \straindata) + \blambda &= \zero, \\
    &\delta\stress:&
    \iweight (\stress - \stressdata) + \bBcoor \nullproj \bmu &= \zero, \\
    &\delta\blambda:&
    \strain - \straincoor &= \zero, \\
    &\delta\bmu:&
    \nullprojT (\bBcoorT \stress - \force) &= \zero, \\
    &\delta\bnu:&
    \cons &= \zero.
  \end{align}
\end{subequations}
Notice that $\strain$ and $\blambda$ can be eliminated by direct substitution
of $\strain = \straincoor$ and $\blambda =-\weight(\strain-\straindata)$.
However, in order to study the overall problem
we are not going to eliminate anything unless strictly necessary.

Now, the second term of the first equation can be evaluated as
\begin{equation}
  \begin{split}
    (\pd{\coor} [\nullprojT (\bBcoorT \stress - \force)])^T \bmu
    &=
    [\bU_2(\coor, \stress)^T \nullproj +
      \bW_1(\coor, \bBcoorT \stress - \force)] \bmu \\
    &=
    \bU_2(\coor, \stress)\nullproj \bmu +
    \bW_2(\coor, \bmu)^T (\bBcoorT \stress - \force),
  \end{split}
\end{equation}
where
\begin{subequations}
  \begin{align}
    \bU_2(\coor, \stress) &:= \pd{\coor}(\bBcoorT \stress) =
    \bU_2(\coor, \stress)^T, \\
    \bW_1(\coor, \ba) &:= \pd{\coor}(\nullprojT \ba), \\
    \bW_2(\coor, \bmu) &:= \pd{\coor}(\nullproj \bmu).
  \end{align}
\end{subequations}
The linearization of the variation of $\lagrangean\sfix$
can be expressed as usual as
\begin{equation}
  \Delta \delta\lagrangean\sfix\argvar
  =
  \delta\bx\sfix \sdot
  \KKT\sfix(\coor, \stress, \blambda, \bmu, \bnu)
  \Delta\bx\sfix
\end{equation}
with the primal-dual NLP variable vector
\begin{equation}
  \bx\sfix
  :=
  (\coor^T, \strain^T, \stress^T, \blambda^T, \bmu^T, \bnu^T)^T.
\end{equation}
The KKT matrix $\KKT\sfix$ is symmetric indefinite. It can be written as
\begin{equation}
  \KKT\sfix(\coor, \stress, \blambda, \bmu, \bnu)
  =
  \begin{bmatrix}
    \KKTcoorcoor & \zero & \KKTstresscoor^T &
    -\bBcoorT & \KKTmucoor^T & \bHcoorT \\
    \zero & \weight & \zero & \bI & \zero & \zero \\
    \KKTstresscoor & \zero & \iweight & \zero & \KKTmustress^T & \zero \\
    -\bBcoor & \bI & \zero & \zero & \zero & \zero \\
    \KKTmucoor & \zero & \KKTmustress & \zero & \zero & \zero \\
    \bHcoor & \zero & \zero & \zero & \zero & \zero
  \end{bmatrix}
\end{equation}
where
\begin{subequations}
  \begin{align}
    \KKTcoorcoor &:=
    -\bU_2(\coor, \blambda) + \bU_2(\coor, \stress) \bW_2(\coor, \bmu) +
    \bW_2(\coor, \bmu)^T \bU_2(\coor, \stress) + \bV(\coor, \bnu), \\
    \KKTstresscoor &:=
    \bU_1(\coor, \nullproj \bmu) + \bBcoor \bW_2(\coor, \bmu), \\
    \KKTmucoor &:=
    \nullprojT \bU_2(\coor, \stress) +
    \bW_1(\coor, \bBcoorT \stress - \force), \\
    \KKTmustress &:= \nullprojT \bBcoorT, \\
    \bU_1(\coor, \ba) &:= \pd{\coor}(\bBcoor\ba), \\
    \bV(\coor, \bnu) &:= \pd{\coor}(\bHcoorT\bnu) = \bV(\coor, \bnu)^T.
  \end{align}
\end{subequations}
Even though explicit expressions for these definitions
depend on the kinematic description adopted,
the format proposed here is generally valid.
As the KKT matrix $\KKT\sfix$ is non-singular,
all local minima of \fixNLP\ are strict minima.
Hence, standard nonlinear optimization methods
that reduce to a full-step Newton iteration in the local convergence area
are highly suited for solving \fixNLP.

The overall DCNLP is often treated by meta-heuristic methods.
However, since it has no useful structure
with respect to the discrete variables
$(\straindata, \stressdata) \in \datacmanifold$,
a mathematically rigorous solution requires enumeration,
that is, finding the minimal value over all measurements
$(\straindata, \stressdata) \in \datacmanifold$
by solving every \fixNLP\ globally.
Alternatively, a well-behaved MINLP formulation
could be solved by rigorous mathematical methods,
but this is highly expensive
already for relatively ``easy'' special cases such as the one discussed in \cite{Kanno2019}.
Therefore we suggest a different approach:
to add suitable structure that enables us
to replace the DCNLP with a single approximating NLP.

\subsection{An approximate nonlinear optimization problem}

An approximate NLP can be stated as
\begin{equation}
  \begin{split}
    \min_{(\strainrec,\stressrec,\coor,\strain,\stress)} \quad &
    \frac12 \|\strain - \strainrec\|^2_{\weight} +
    \frac12 \|\stress - \stressrec\|^2_{\iweight} \\
    \text{subject to} \quad
    &\strain - \straincoor = \zero, \\
    &\nullprojT (\bBcoorT \stress - \force) = \zero, \\
    &\cons =\zero, \\
    &\bg(\strainrec,\stressrec) = \zero.
  \end{split}
\end{equation}
The new strain and stress variables $\strainrec$ and $\stressrec$
are parameters that describe an underlying constitutive manifold
in its reconstructed version (an approximation),
which is implicitly defined as
\begin{equation}
  \reccmanifold
  :=
  \{(\strainrec, \stressrec) \in \reals^{\nstrain} \mid
  \bg(\strainrec,\stressrec) = \zero \in \reals^{\nstrain}\}
\end{equation}
and that satisfies
\begin{equation}
  \|\bg(\straindata, \stressdata)\| \le \varepsilon
  \quad \forall \, (\straindata, \stressdata) \in \datacmanifold
\end{equation}
for some accuracy tolerance $\varepsilon > 0$.
Additionally, physical consistency requires that
$\bg(\strainrec, \zero) = \zero$ implies $\strainrec = \zero$ and
$\bg(\zero, \stressrec) = \zero$ implies $\stressrec = \zero$.
The idea here is to replace the measurement data set
by enforcing the state to belong to the reconstructed constitutive manifold
that has a precise mathematical structure
and that is derived from the same data set.
The underlying assumption is, of course, that such a constitutive manifold exists
and that we can reconstruct a (smooth) implicit representation $\bg$
whose linearization takes the simple form
\begin{equation}
  \pd{\strainrec} \bg(\strainrec, \stressrec) \Delta\strainrec +
  \pd{\stressrec} \bg(\strainrec, \stressrec) \Delta\stressrec = \zero,
\end{equation}
which can be interpreted as a sort of hidden constraint. For materials with symmetric properties,
a further condition can be stated as
\begin{equation}
  %[\pd{\stressrec} \bg(\strainrec, \stressrec)]^{-1}
  %[\pd{\strainrec} \bg(\strainrec, \stressrec)] -
  %[\pd{\strainrec} \bg(\strainrec, \stressrec)]^T
  %[\pd{\stressrec} \bg(\strainrec, \stressrec)]^{-T}
  %=
  %\zero,
  \mathrm{skew}\big([\pd{\stressrec} \bg(\strainrec, \stressrec)]^{-1}[\pd{\strainrec} \bg(\strainrec, \stressrec)]\big) = \zero
  \quad \textrm{or} \quad
  \mathrm{skew}\big([\pd{\strainrec} \bg(\strainrec, \stressrec)]^{-1}[\pd{\stressrec} \bg(\strainrec, \stressrec)]\big) = \zero,
\end{equation}
requiring the regularity of $\pd{\strainrec} \bg(\strainrec, \stressrec)$
and $\pd{\stressrec} \bg(\strainrec, \stressrec)$.
The reconstructed constitutive manifold will enormously
facilitate the task of the \dd~solver,
avoiding the cost of solving a DCNLP,
either by enumeration
or by heuristic or meta-heuristic methods. The latter can in general only provide approximate solutions
that strongly depend on the initial guess and whose convergence properties
are inferior when compared to gradient-based methods.

A constitutive manifold is said to be thermomechanically consistent
if it is derived from a hyperelastic energy function $\Psi$
such that the following functional structure holds
\cite{Ibanez2017,Crespo2017}:
\begin{equation}
  \bg(\strainrec, \stressrec)
  =
  \stressrec - \pd{\strainrec} \Psi(\strainrec)
  =
  \zero.
\end{equation}
The reconstruction of the hyperelastic energy function $\Psi$
is very attractive due to three main reasons:
\textit{i}) it ensures the thermomechanical consistency
and then, all symmetries are retained;
\textit{ii}) reconstructing a single scalar function represents a smaller computational expense
than reconstructing all components of the elastic tensor; and,
\textit{iii}) no reformulation of the finite element technology is necessary.
However, in the case of new composite materials
or meta materials that exhibit non-convex responses,
the reconstruction of the energy function may not be very convenient.
More importantly, in some cases the formulation
of an energy function may not even be possible.
Thus, we adopt the constitutive manifold $\reccmanifold$
as introduced previously without assuming any special functional structure
of the constitutive constraint $\bg$.
Further specializations are possible and should be instantiated
for specific applications of the proposed formulation.

The Lagrangian function of the approximate NLP is then given by
\begin{equation}
  \begin{split}
    \lagrangean\sone\argvarrec
    &=
    \lagrangean\sfix(\bx\sfix; \strainrec, \stressrec) +
    \bxi \sdot \bg(\strainrec, \stressrec) \\
    &=
    \frac12 \|\strain - \strainrec\|^2_{\weight} +
    \frac12 \|\stress - \stressrec\|^2_{\iweight} \\
    &+ \blambda \sdot (\strain - \straincoor) \\
    &+ (\nullproj \bmu) \sdot (\bBcoorT \stress - \force) \\
    &+ \bnu \sdot \cons \\
    &+ \bxi \sdot \bg(\strainrec, \stressrec),
  \end{split}
\end{equation}
where $\bxi \in \reals^{\nstrain}$ are Lagrange multipliers
that correspond to the enforcement of the strain and stress states
to remain on the constitutive manifold.

Once again, to find the first-order optimality conditions,
the variation of $\lagrangean\sone$ is calculated as:
\begin{equation}
  \begin{split}
    \delta\lagrangean\sone\argvarrec
    &= \delta\strainrec \sdot
    (\weight (\strainrec - \strain) +
    [\pd{\strainrec} \bg(\strainrec, \stressrec)]^T \bxi) \\
    &+ \delta\stressrec \sdot
    (\iweight (\stressrec - \stress) +
    [\pd{\stressrec} \bg(\strainrec, \stressrec)]^T \bxi) \\
    &+ \delta\coor \sdot (-\bBcoorT \blambda +
          (\pd{\coor} [\nullprojT (\bBcoorT\stress - \force)])^T \bmu +
          \bHcoorT \bnu) \\
    &+ \delta\strain \sdot (\weight (\strain - \strainrec) + \blambda) \\
    &+ \delta\stress \sdot
    (\iweight (\stress - \stressrec) + \bBcoor \nullproj \bmu) \\
    &+ \delta\blambda \sdot (\strain - \straincoor) \\
    &+ \delta\bmu \sdot \nullprojT (\bBcoorT \stress - \force) \\
    &+ \delta\bnu \sdot \cons \\
    &+ \delta\bxi \sdot \bg(\strainrec, \stressrec).
  \end{split}
\end{equation}
Setting this to zero for any choice of the varied quantities,
we obtain the following KKT conditions:
\begin{subequations}
  \begin{align}
    &\delta\strainrec: &
    -\weight (\strain - \strainrec) +
    [\pd{\strainrec} \bg(\strainrec, \stressrec)]^T \bxi &= \zero, \\
    &\delta\stressrec: &
    -\iweight (\stress - \stressrec) +
    [\pd{\stressrec} \bg(\strainrec, \stressrec)]^T \bxi &= \zero, \\
    &\delta\coor: &
    -\bBcoorT \blambda +
    (\pd{\coor} [\nullprojT (\bBcoorT \stress - \force)])^T \bmu +
    \bHcoorT \bnu &= \zero, \\
    &\delta\strain: &
    \weight (\strain - \strainrec) + \blambda &= \zero, \\
    &\delta\stress: &
    \iweight (\stress - \stressrec) + \bBcoor \nullproj \bmu &= \zero, \\
    &\delta\blambda: & \strain - \straincoor &= \zero, \\
    &\delta\bmu: & \nullprojT (\bBcoorT \stress - \force) &= \zero, \\
    &\delta\bnu: & \cons &= \zero, \\
    &\delta\bxi: & \bg(\strainrec, \stressrec) &= \zero.
  \end{align}
\end{subequations}
The linearization of the variation of $\lagrangean\sone$
can be expressed as
\begin{equation}
  \Delta\delta\lagrangean\sone\argvarrec
  =
  \delta\bx\sone \sdot
  \KKT\sone
  (\strainrec, \stressrec, \coor, \stress, \blambda, \bmu, \bnu, \bxi)
  \Delta\bx\sone
\end{equation}
with
\begin{equation}
  \bx\sone := (\strainrec^T, \stressrec^T, \bx\sfix^T, \bxi^T)^T.
\end{equation}
The explicit form of the KKT matrix $\KKT\sone$ is
\begin{equation}
  \KKT\sone(\strainrec, \stressrec, \coor, \stress, \blambda, \bmu, \bnu, \bxi)
  =
  \begin{bmatrix}
    \KKT_{\strainrec \strainrec}(\strainrec, \stressrec, \bxi) &
    \KKT_{\stressrec \strainrec}(\strainrec, \stressrec, \bxi)^T &
    \KKT_{\bx \strainrec}^T &
    [\pd{\strainrec} \bg(\strainrec,\stressrec)]^T \\
    \KKT_{\stressrec \strainrec}(\strainrec, \stressrec, \bxi) &
    \KKT_{\stressrec \stressrec}(\strainrec, \stressrec, \bxi) &
    \KKT_{\bx \stressrec}^T &
    [\pd{\stressrec} \bg(\strainrec,\stressrec)]^T \\
    \KKT_{\bx \strainrec} &
    \KKT_{\bx \stressrec} &
    \KKT\sfix(\coor, \stress, \blambda, \bmu, \bnu) &
    \zero \\
    \pd{\strainrec} \bg(\strainrec, \stressrec) &
    \pd{\stressrec} \bg(\strainrec, \stressrec) &
    \zero &
    \zero \\
  \end{bmatrix}
\end{equation}
where
\begin{subequations}
  \begin{align}
    \KKT_{\strainrec \strainrec}(\strainrec, \stressrec, \bxi)
    &:=
    \weight +
    \pd{\strainrec} ([\pd{\strainrec} \bg(\strainrec, \stressrec)]^T \bxi) =
    \KKT_{\strainrec \strainrec}(\strainrec, \stressrec, \bxi)^T, \\
    \KKT_{\stressrec \strainrec}(\strainrec, \stressrec, \bxi)
    &:=
    \pd{\stressrec} ([\pd{\strainrec} \bg(\strainrec, \stressrec)]^T \bxi), \\
    \KKT_{\stressrec \stressrec}(\strainrec, \stressrec, \bxi)
    &:=
    \iweight +
    \pd{\stressrec} ([\pd{\stressrec} \bg(\strainrec, \stressrec)]^T \bxi) =
    \KKT_{\stressrec \stressrec}(\strainrec, \stressrec, \bxi)^T, \\
    \KKT_{\bx \strainrec}
    &:=
    \begin{bmatrix} \zero & -\weight & \zero & \zero & \zero \end{bmatrix}^T, \\
    \KKT_{\bx \stressrec}
    &:=
    \begin{bmatrix} \zero & \zero &-\iweight & \zero & \zero \end{bmatrix}^T.
  \end{align}
\end{subequations}
Again the KKT matrix $\KKT\sone$ is non-singular,
hence all local minima are strict and the approximate NLP can be solved robustly.
In the following section, we provide an example that demonstrates
how to derive the concrete approximate NLP
for a given finite element formulation.

%% file: 4structure.tex
\section{Application of the proposed approach}

In this section, we describe two structural models
that are reformulated within the proposed setting of \DDCM.
The first model is a \dd~truss element that serves as a starting point;
this element type was already successfully investigated in
\cite{Kirchdoerfer2016,Nguyen2018}.
In contrast to the approaches available in the literature,
we apply the framework developed in the previous section to the truss element,
unveiling details of its global format that to the best of our knowledge
have not yet been published elsewhere.
Therefore we call this an ``illustrative example''.
The second model, a main innovation of the present work,
is a \dd~geometrically exact beam that is given in a
frame-invariant path-independent finite element formulation.
This model relies on a kinematically constrained approach,
where the orientation of the cross section is described by means of
three vectors that are constrained to be mutually orthonormal.
Both examples have favorable mathematical structures
that are exploited to derive the required
finite element machinery analytically.

\subsection{\Dd~truss element (illustrative example)}

The position of any point belonging to a truss element can be written as
\begin{equation}
  \bphi(x) = \left(1 - \frac{x}{L}\right) \bphi_0 + \frac{x}{L} \bphi_L
  \in \reals^3,
\end{equation}
with $\bphi_0$ and $\bphi_L$ in $\reals^3$ being the positions of both ends,
in which $L$ is the reference length and $x$ is a spatial variable.
By defining
\begin{equation}
  \coor =
  \begin{pmatrix}
    \bphi_0 \\
    \bphi_L
  \end{pmatrix},
\end{equation}
we can compute the axial Green-Lagrange strain as
\begin{equation}
  e\ofcoor = \frac12 \coor^T \! \bm{E} \coor - e_\refe,
\end{equation}
where $\bE \in \reals^{6 \times 6}$ is a symmetric matrix
and the subscript ``$\refe$'' indicates the stress-free configuration.

The operator $\bBcoor = \pd\coor e\ofcoor$,
which relates the variation of displacement-based strains
with the variation of the kinematic fields through the relation
$\delta e = \bBcoor \delta\coor$, has the explicit form
\begin{equation}
  \bBcoor = \coor^T \! \bE,
\end{equation}
and its derivative with respect to $\coor$ is
\begin{equation}
  \pd{\coor} \bBcoor = \bE.
\end{equation}
For simplicity, we consider constraints that are linear in $\coor$,
\begin{equation}
  \cons = \bm{H} \coor - \bm{h}_\refe,
\end{equation}
where $\bh_\refe$ is a reference offset vector
and $\bH$ is the constant Jacobian matrix.
The null-space matrix $\bm{N}$ is also constant and can be determined
by inspection once the constraints are specified.
For this case the cost function to be minimized adopts the very simple form
\begin{equation}
  \cost^{\mathrm{truss}}(\tilde{e}, \tilde{s}, e, s)
  =
  \frac12 c (e - \tilde{e})^2 + \frac12 c^{-1} (s - \tilde{s})^2,
\end{equation}
where the scalar $c$ is just a weight factor.
Finally, by putting all pieces together,
we obtain the KKT conditions of \fixNLP:
\begin{subequations}
  \begin{align}
    &\delta\coor:&
    -\lambda \bm{E} \coor + s \bm{E} \bm{N} \bmu + \bm{H}^T \! \bnu &= \zero, \\
    &\delta e:&
    c (e - \tilde{e}) + \lambda &= 0, \\
    &\delta s:&
    c^{-1} (s - \tilde{s}) + \bm{q}^T \! \bm{E} \bm{N} \bmu &= 0, \\
    &\delta\lambda:&
    -\frac12 \coor^T \! \bm{E} \coor + e + e_\refe &= 0, \\
    &\delta\bmu:&
    s \bm{N}^T \! \bm{E} \coor - \bm{N}^T \! \force &= \zero, \\
    &\delta\bnu:&
    \bm{H} \coor - \bm{h}_\refe &= \zero.
  \end{align}
\end{subequations}
Linearization yields the KKT matrix
(which does not depend on $e$ and $\bnu$ here),
\begin{equation}
  \KKT\sfix^{\mathrm{truss}}(\coor, s, \lambda, \bmu)
  =
  \begin{bmatrix}
    -\lambda \bm{E} & \zero & \bm{E} \bm{N} \bmu &
    -\bm{E} \coor & s \bm{E} \bm{N} & \bm{H}^T \\
    \zero & c & 0 & 1 & \zero & \zero \\
    \bmu^T \! \bm{N}^T \! \bm{E} & 0 & c^{-1} & 0 &
    \coor^T \! \bm{E} \bm{N} & \zero \\
    -\coor^T \! \bm{E} & 1 & 0 & 0 & \zero & \zero \\
    s \bm{N}^T \! \bm{E} & \zero & \bm{N}^T \! \bm{E} \coor &
    \zero & \zero & \zero \\
    \bm{H} & \zero & \zero & \zero & \zero & \zero
  \end{bmatrix}
  .
\end{equation}
Formulating the approximate NLP is straightforward
and therefore omitted here.
However, when looking at the DCNLP approach, we can identify some drawbacks
that confirm the strength of our newly proposed approach.

Let us consider a single truss element of unit initial length.
One end is rigidly fixed and the other end is free to move
in the longitudinal direction, \ie we have
$\bH = [\bI_{5\times5}, \zero_{5\times1}]$ and
$\bm h_\refe = \zero_{5\times1}$.
At the free end we consider a longitudinal force with magnitude $20$.
By considering a unit stiffness,
$\bE = \frac12[\bI_{3\times3},-\bI_{3\times3};-\bI_{3\times3},\bI_{3\times3}]$
and $e_\refe = 1$, the equilibrium can be explicitly written as
\begin{equation}
  f^{\textrm{int}} - f^{\textrm{ext}} = 0 \implies \frac12 (q^3 - q) - 20 = 0,
\end{equation}
where $q$ is the single absolute coordinate
that describes the motion of the free end.
This equation is easily solved; it possesses three roots.
One root is purely real, $q = 3.51739351$,
and the two remaining ones are complex conjugates,
which lack physical meaning.
The equilibrium state corresponds to the strain-stress pair
$(5.68602856, 5.68602856)$.

\begin{figure}[t]
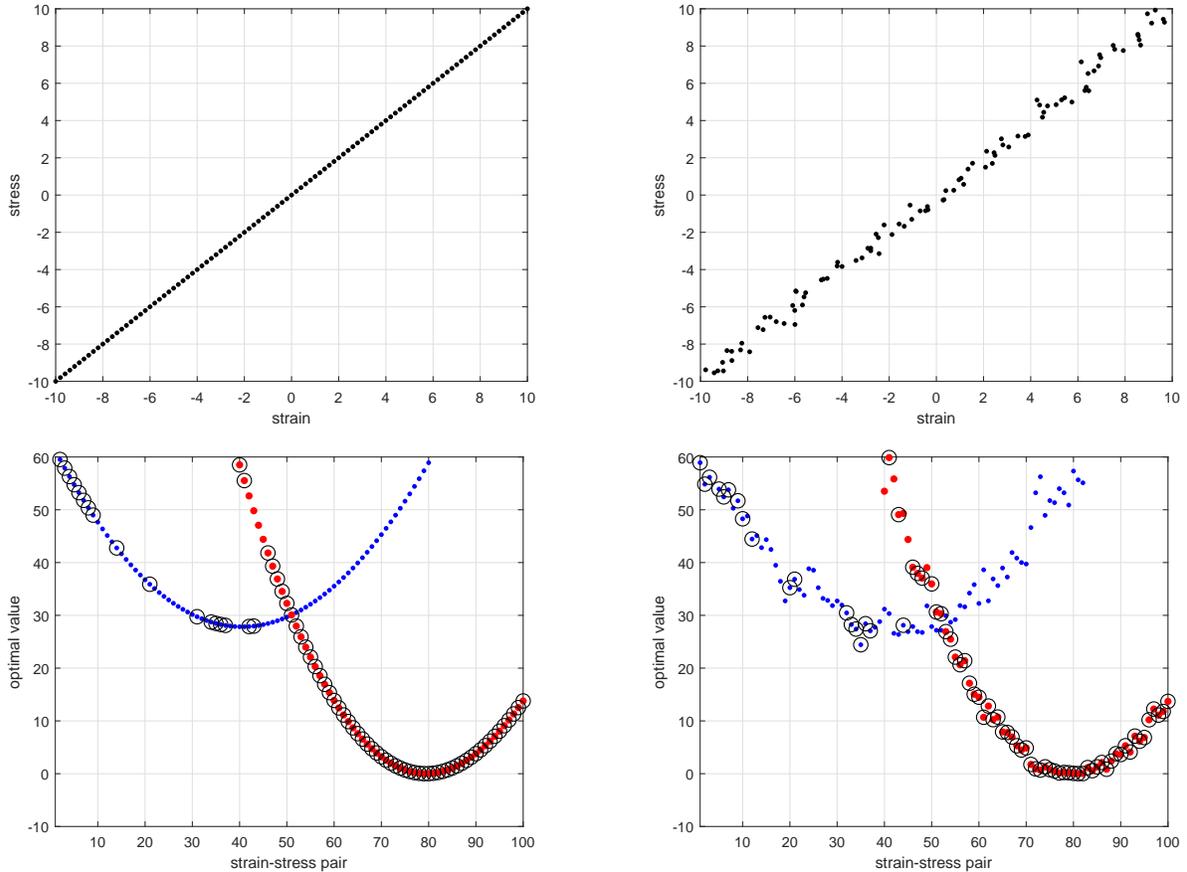

  \centering
  \dfig[45]1\hfill\sfig[45]1\\[2ex]
  \dfig[45]2\hfill\sfig[45]2
  \caption{Synthetic data set and optimal values
    of the cost function vs.\ synthetic data;
    exact (left), 10\% noise (right).
    Small black dots indicate points in the phase space;
    black circles indicate cold start,
    small blue dots indicate ascending warm start,
    and large red dots indicate descending warm start.}
  \label{fig:tlnlptruss}
\end{figure}

\begin{figure}[tp]
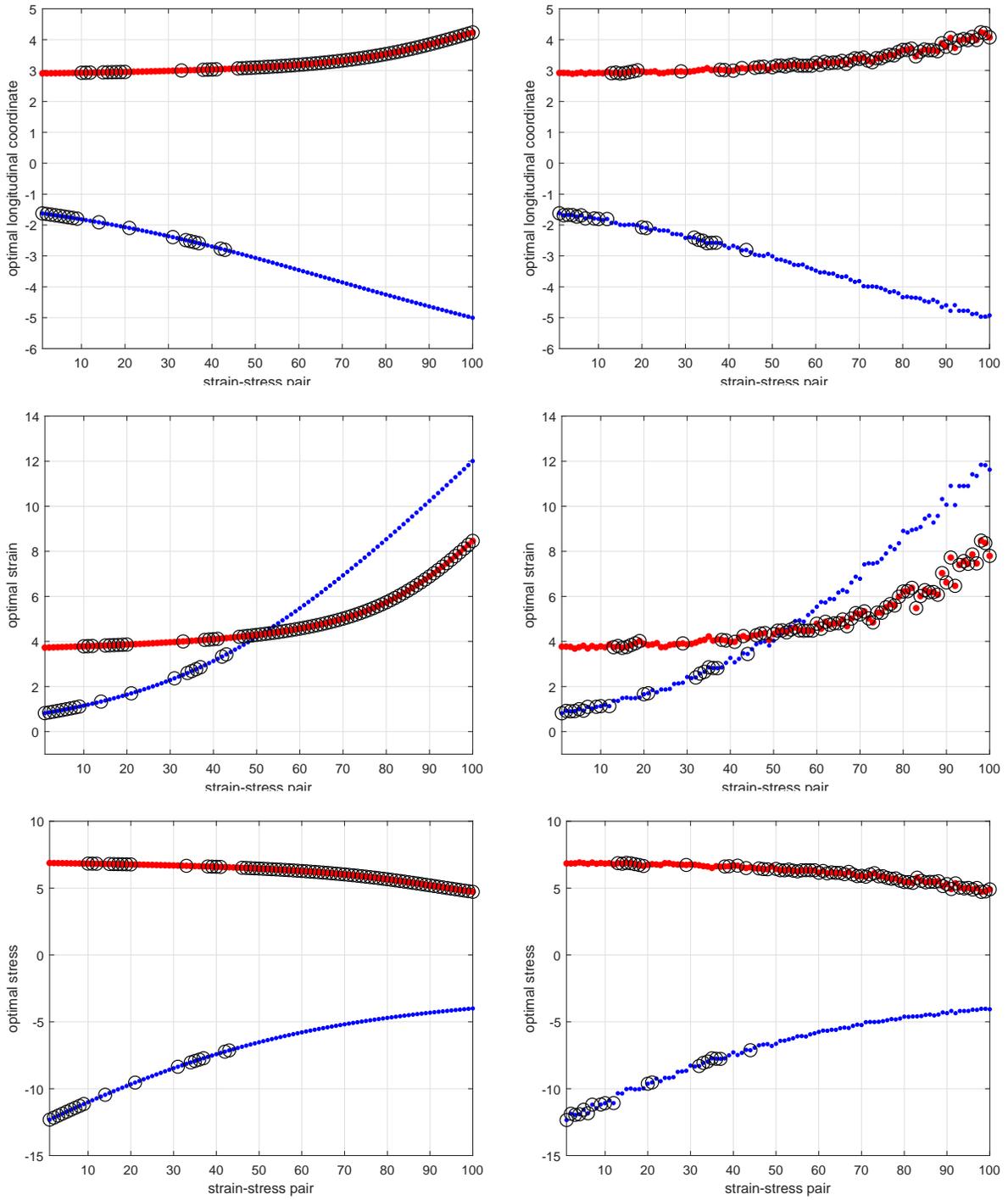

  \centering
  \dfig3\hfill\sfig3\\[2ex]
  \dfig4\hfill\sfig4\\[2ex]
  \dfig5\hfill\sfig5
  \caption{Optimal states vs.\
    exact synthetic data set (left) and with 10\% noise (right).
    Black circles indicate cold start,
    small blue dots indicate ascending warm start,
    and large red dots indicate descending warm start.}
  \label{fig:tr-opt}
\end{figure}

Now, we explore how the DCNLP approach works
by solving every \fixNLP\ in turn.
To this end, we generate two synthetic data sets:
one that lies exactly on the manifold ($e = s$)
and one with 10\% noise added,
each consisting of $101$ pairs that are stored in ascending order,
see Figure~\ref{fig:tlnlptruss} (top).
Then, we solve \fixNLP\ with $c = 1$
for each strain-stress pair in the data set
by a local SQP method
(\ie a full-step Newton iteration on the KKT conditions)
starting from the stress-free state ($q = 1$)
with zero Lagrange multipliers (cold start, indicated by black circles)
and with a relative error-based tolerance of $10^{-10}$,
\ie we iterate until $\|\Delta\bx\| \le 10^{-10} \|\bx\|$.
In addition, we solve each \fixNLP\ warm-started from
the previous solution in ascending order (indicated by small blue dots)
and in descending order of data pairs (indicated by large red dots),
see Figure~\ref{fig:tlnlptruss} (bottom).
We observe that every \fixNLP\ has two distinct local minimizers.
It can be shown that we actually have two manifolds of local minimizers
when considering \fixNLP\ as a parametric NLP whose parameters
$(\straindata, \stressdata)$ vary in the manifold $e = s$.
The cold-started Newton iteration does not converge in all cases,
and if it converges one cannot predict to which minimizer.
In contrast, both warm-started Newton iterations always converge
to the ``same'' minimizer because the data points are sufficiently close.
When cold-started, the solver usually takes six to eight iterations,
but sometimes more than ten iterations in the vicinity of the
``cross-over'' between the two minimizers.
The iteration count until convergence reduces to four with the warm start.
The global DCNLP minimum corresponds to the pair number $79$,
$(5.60000000, 5.60000000)$,
which is the strain-stress pair that is closest to the analytic solution.
\ifcase1
% Show Lagrange multipliers in plots.
Figures \ref{fig:tr-opt-00} and~\ref{fig:tr-opt-10} reveal
\begin{figure}[tp]
  \centering
  \dfig3\hfill\dfig4\\[2ex]
  \dfig5\hfill\dfig6\\[2ex]
  \dfig7\hfill\dfig8
  \caption{Optimal states and multipliers
    vs.\ synthetic data set without noise.}
  \label{fig:tr-opt-00}
\end{figure}
\begin{figure}[tp]
  \centering
  \sfig3\hfill\sfig4\\[2ex]
  \sfig5\hfill\sfig6\\[2ex]
  \sfig7\hfill\sfig8
  \caption{Optimal states and multipliers
    vs.\ synthetic data set with 10\% noise.}
  \label{fig:tr-opt-10}
\end{figure}
\else
% Omit Lagrange multipliers.
Figure~\ref{fig:tr-opt} reveals
\fi
that the longitudinal coordinate $q$ and the stress $s$
are always positive for one minimizer (large red dots)
and negative for the other one (small blue dots);
hence the truss element is under traction in the first case, as intended,
and it is under (physically impossible) compression in the second case.
It is well-known that the mechanical behavior provided by the
Saint Venant--Kirchhoff material model for the total compression case
disagrees with experimental observations.

This very simple example shows the difficulties that ought to be faced
when using the DCNLP approach with unprocessed data.
Each \fixNLP\ can have several local minimizers.
Some of these minimizers, including the global one,
may be physically impossible.
The full-step Newton method used above may not converge in all cases,
but then one can use more sophisticated local NLP solvers
that employ a line search or trust region
and that are guaranteed to converge
as long as numerical roundoff errors are sufficiently small.
Although one would expect the global DCNLP minimizer
to be the physically correct one,
% \MCS[I believe it is unique?] Yes, in all cases that we consider.
which turns out to be true in our example,
it is unclear whether this property must always hold.
% \MCS[It is probably true in all literature examples?]
% Unclear: no rigorous methods used.
\ifcase1% Irrelevant here.
In any case, sufficiently poor measurements
can certainly destroy this property.
\fi
Provided that the global DCNLP minimizer is the physically correct one,
solving DCNLP rigorously is very hard (actually NP-hard), as mentioned.
The enumeration approach explored above requires the global solution
of at least one \fixNLP\ (the expensive part).
Then, that parametric minimizer can simply be tracked
by inexpensive warm-started solves
through sequences of sufficiently close data pairs.
\ifcase1% Elaborate (unclear):
It is perhaps quite likely that a first global minimizer can be found
by solving several \fixNLP\ locally from the same starting point.
\fi
Alternatively, a well-behaved MINLP formulation would be tractable
by rigorous solvers, which is known to be similarly expensive.
These rigorous solvers, however, have the essential advantage that
intermediate near-optimal solutions come with a quality measure.
In recent work \cite{Kirchdoerfer2016, Kirchdoerfer2018, Nguyen2018},
meta-heuristic optimizers, like simulated annealing,
are employed to solve the DCNLP.
It is well-known that these methods
approximate the solution without any quality measure,
have poor convergence properties when compared to gradient-based optimizers,
and the solution found depends strongly on the initial guess.
On the other hand, a constitutive manifold behind the data
always exists by physical reasons,
and we %strongly
suggest to reconstruct and use it explicitly.
Outside the optimization context, Iba\~nez \etal \cite{Ibanez2019} suggest to take advantage
of the acquired knowledge on the constitutive characterization of materials and propose,
by means of machine learning techniques, to develop corrections to popular material models.
After the initial expensive reconstruction step,
it reduces the influence of measurement errors on the minimizers,
it guarantees the possibility of warm starts,
it allows non-standard elasticity models,
it enables gradient-based solvers, and
it reduces the computational effort.
In addition, global minimizers are easily detected
by their zero residual: they lie on the manifold.
Together with the aspects discussed above,
these are strong arguments in favor of our approximate NLP approach.
Finally, any rigorous optimization approach offers
the advantage of adding inequality constraints.
Given sufficient physical insight, this can be used
to ``cut off'' some or all non-physical minimizers,
such as the one with certain negative states in the example.
This can drastically reduce the solution effort.

\subsection{\Dd~geometrically exact beam element}

\begin{figure}[tp]
  \centering
  \includegraphics[trim=30 60 30 80,clip,width=0.85\textwidth]{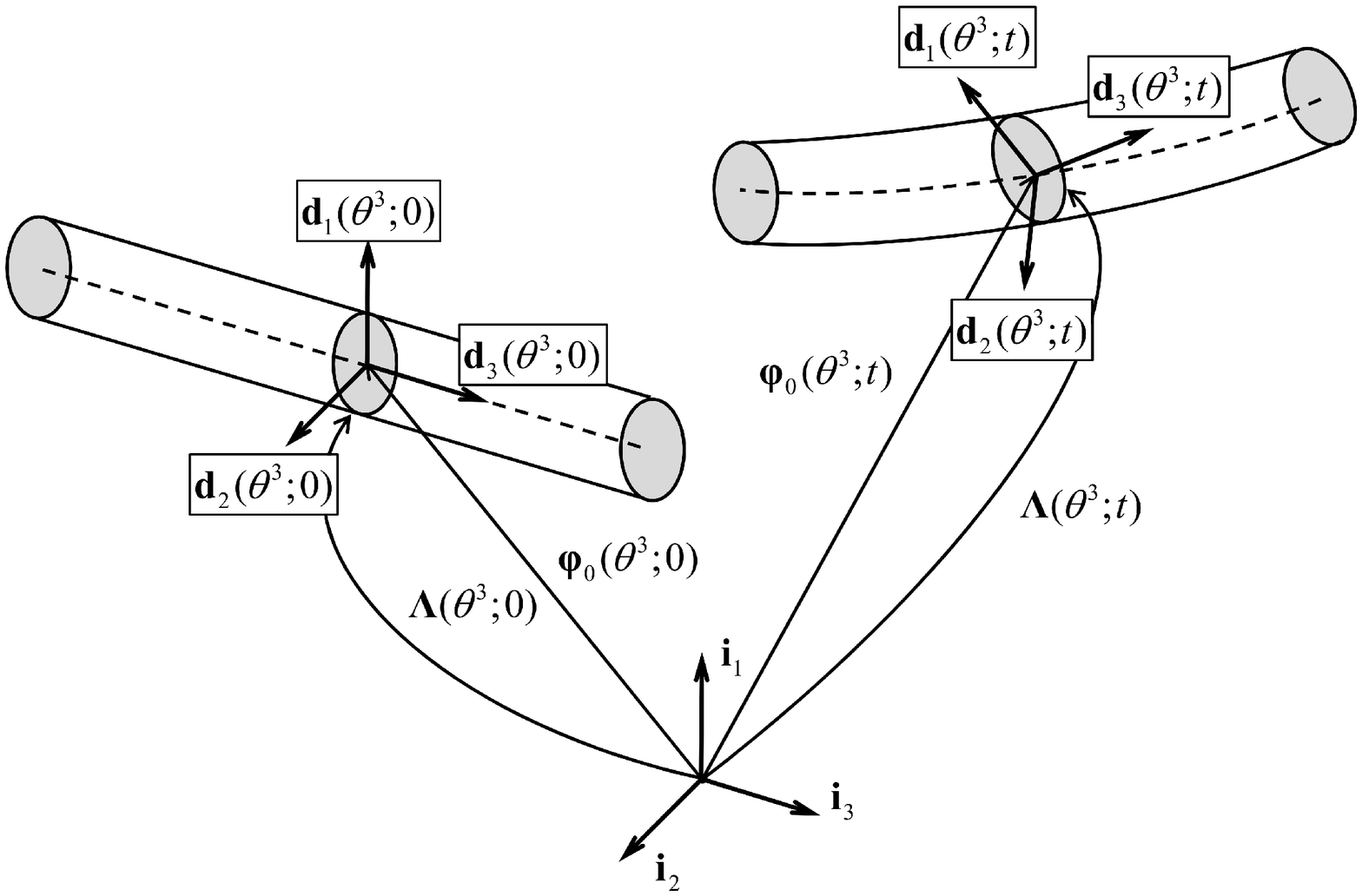}
  \caption{The geometrically exact beam:
    evolution among configurations through the regular motion
    $\bphi(\rdot; t) \circ \bphi(\rdot;0)^{-1}$.}
  \label{fig:beam}
\end{figure}

The position of any point belonging to the beam
shown in Figure~\ref{fig:beam} can be written as
\begin{equation}
  \bphi(\btheta)
  =
  \coorz(\theta^3) + \theta^1 \coori(\theta^3) + \theta^2 \coorj(\theta^3)
  \in \reals^3,
\end{equation}
in which $\coorz \in \reals^3$ is the position vector of the beam axis
and $\coori \in S^2$, $\coorj \in S^2$ together with $\coork \in S^2$
are three mutually orthonormal directors.

The directors can be described by means of the unit sphere,
which is a nonlinear, smooth, compact, two-dimensional manifold
that can be embedded in $\reals^3$ as
\begin{equation}
  \label{eq-stwo}
  S^{2} := \{\bm{d} \in \reals^3 \mid \bm{d} \sdot \bm{d} = 1\}.
\end{equation}
Special attention must be paid to the fact that this
manifold possesses no special algebraic structure,
specifically group-like structure \cite{Eisenberg1979}.
On that basis, the rotation tensor for the cross section is simply obtained as
$\rotation = \coori\otimes\uii+\coorj\otimes\ujj+\coork\otimes\ukk \in SO(3)$,
in which $\{\uii, \ujj, \ukk\}$ is the dual basis
of the ambient space $\euclideans^3$
($\reals^3$ with the standard Euclidean structure),
\ie the basis of the space of row vectors.
The group of rotations is a nonlinear, smooth, compact,
three-dimensional manifold defined as
\begin{equation}
  SO(3) :=
  \{
  \rotation \in \reals^{3 \times 3} \mid
  \rotation^T \rotation = \bI,
  \det\rotation = 1
  \}.
\end{equation}
In contrast to the unit sphere,
this manifold does possess a group-like structure
when considered with the tensor multiplication operation,
hence it is a Lie group.

%Its Lie algebra is  the set $so(3)$ of skew-symmetric $3\times3$ matrices.  An isomorphism exists between vectors in $\reals^3$ and $so(3)$ defined as $\hat{\cdot}:\reals^3\to so(3)$ such that for all $\mbs{w},\mbs{a}\in\reals^3$, the tensor $\hat{\mbs{w}}\in so(3)$ satisfies $\hat{\mbs{w}}\mbs{a} = \mbs{w}\times \mbs{a}$.  The vector $\mbs{w}$ is referred to as the axial vector of the skew-symmetric tensor $\hat{\mbs{w}}$ and we also write $\skew[\mbs{w}] = \hat{\mbs{w}}$.  The exponential map $\exp:so(3)\to SO(3)$ is a surjective application.
The set of parameters $\btheta = (\theta^1, \theta^2, \theta^3)$
is chosen in a way that the vector
$\bar{\btheta} = \theta^1 \coori + \theta^2 \coorj$
completely describes the cross section.
%, which intersects the beam axis and $\bphi=\bphi(\btheta)$ is the given parametrization rule in space.
In the context of geometrically exact beams,
the doubly-covariant displacement-based Green-Lagrange strain tensor, $\gl(\bphi)$,
can be simplified by eliminating quadratic strains.
Thus, its components are approximated as
\begin{equation}
  E_{ij}%(\btheta)
  \approx
  \mathrm{symm}(\delta_{i3}\delta_{jk}((\gamma^k%(\theta^3)
  -\gamma^k_{\refe}%(\theta^3)
  )-\epsilon^k_{lm}\bar{\theta}^l(\omega^m%(\theta^3)
  -\omega^m_{\refe}%(\theta^3)
  ))),
\end{equation}
where $\mathrm{symm}(\rdot)$ stands for the symmetrization
of the tensor considered.
The subscript ``$\refe$'' indicates the stress-free configuration,
$\delta_{ij}$ denotes the Kronecker delta,
and $\epsilon^i_{jk}$ is the alternating symbol
that appears in the computation of the cross product
in three-dimension Euclidean space.
From now on, we set $\theta^3 = \arclen$ to indicate every reference
related to the arc length of the beam.
The scalars $\gamma^i$ are the components of a first deformation vector
defined as
\begin{equation}
  \label{eq:beam-gamma}
  \bgamma%\ofcoor
  =
  \begin{pmatrix}
    \coori \sdot \pd{\arclen} \coorz \\
    \coorj \sdot \pd{\arclen} \coorz \\
    \coork \sdot \pd{\arclen} \coorz
  \end{pmatrix}
  .
\end{equation}
For shear refer to first and second components,
and for elongation refer to the third one.
The scalars $\omega^i$ are the components of a second deformation vector
defined as
\begin{equation}
  \label{eq:beam-omega}
  \bomega%\ofcoor
  =
  \frac12
  \begin{pmatrix}
    \coork \sdot \pd{\arclen} \coorj - \coorj \sdot \pd{\arclen} \coork \\
    \coori \sdot \pd{\arclen} \coork - \coork \sdot \pd{\arclen} \coori \\
    \coorj \sdot \pd{\arclen} \coori - \coori \sdot \pd{\arclen} \coorj
  \end{pmatrix}
  .
\end{equation}
For bending refer to first and second components,
and for torsion refer to the third one.
For the sake of compactness, let us introduce the vector
containing all kinematic fields,
\begin{equation}
  \label{eq:beam-kfields}
  \coor(\arclen)
  =
  (\coorz(\arclen)^T, \coori(\arclen)^T, \coorj(\arclen)^T, \coork(\arclen)^T)^T,
\end{equation}
and the vector that gathers the two strain measures
obtained from the kinematic field,
\begin{equation}
  \label{eq:beam-strain}
  \straincoor =
  \begin{pmatrix}
    \bgamma\ofcoor - \bgamma_{\refe} \\
    \bomega\ofcoor - \bomega_{\refe}
  \end{pmatrix}
  .
\end{equation}
Additionally, we need to introduce the vector containing
all generalized ``strain'' fields that is going to be tied
by the compatibility equation,
\begin{equation}
  \strain =
  \begin{pmatrix}
    \bgamma - \bgamma_{\refe} \\
    \bomega - \bomega_{\refe}
  \end{pmatrix}
  ,
\end{equation}
and the vector containing the two generalized ``stress'' fields,
\begin{equation}
  \stress =
  \begin{pmatrix}
    \bforce \\
    \bmoment
  \end{pmatrix}
  ,
\end{equation}
which contains the cross sectional force and moment resultants,
\ie three force components and three moment components.

\ifcase1
The operator that relates the variation of the displacement-based strains
to the variation of the kinematic fields through the relation
$\delta \straincoor = \bBcoor \delta \coor$ has the explicit form
\or
Differentiating the components \eqref{eq:beam-gamma}
and \eqref{eq:beam-omega} of \eqref{eq:beam-strain}
with respect to $\coor$ from \eqref{eq:beam-kfields}
yields the operator $\Boper(\coor) = \pd\coor \straincoor$
in its explicit form
\fi
\begin{equation}
  % MCS: Here we need braces to avoid strange typesetting of some '^T'!?
  \Boper(\coor) = \frac12
  \begin{bmatrix}
    2 \coori^T \pdarclen & 2 \pd{\arclen} {\coorz^T} & \zero & \zero \\
    2 \coorj^T \pdarclen & \zero & 2 \pd{\arclen} {\coorz^T} & \zero \\
    2 \coork^T \pdarclen & \zero & \zero & 2 \pd{\arclen} {\coorz^T} \\
    \zero & \zero & \coork^T \pdarclen - \pd{\arclen} {\coork^T} &
    \pd{\arclen} \coorj^T - \coorj^T \pdarclen \\
    \zero & \pd{\arclen} \coork^T - \coork^T \pdarclen &
    \zero & \coori^T \pdarclen - \pd{\arclen} {\coori^T} \\
    \zero & \coorj^T \pdarclen - \pd{\arclen} {\coorj^T} &
    \pd{\arclen} \coori^T - \coori^T \pdarclen & \zero
  \end{bmatrix}
  .
\end{equation}
It is linear with respect to kinematic fields,
a fact that is very convenient
regarding the computation of high-order derivatives.
For assembling the KKT matrix $\KKT$,
we compute $\Uoper_1(\ba) = \pd\coor (\Boper(\coor) \ba)$ as
\begin{equation}
  \Uoper_1(\ba) = \bUbar_1(\ba) + \bUbre_1(\ba) \Doper,
\end{equation}
where
\begin{math}
  \ba(\arclen)
  =
  (\ba_0(\arclen)^T, \ba_1(\arclen)^T, \ba_2(\arclen)^T, \ba_3(\arclen)^T)^T
  \in \reals^{12}
\end{math}
is a vector with $\ba_i(\arclen) \in \euclideans^3$ for $i = 0, \dots, 3$.
We obtain
\begin{align}
  \bUbar_1(\ba) &= \frac12
  \begin{bmatrix}
    \zero & 2 \pdarclen \ba_0^T &               \zero &               \zero \\
    \zero &               \zero & 2 \pdarclen \ba_0^T &               \zero \\
    \zero &               \zero &               \zero & 2 \pdarclen \ba_0^T \\
    \zero &               \zero &  -\pdarclen \ba_3^T &   \pdarclen \ba_2^T \\
    \zero &   \pdarclen \ba_3^T &               \zero &  -\pdarclen \ba_1^T \\
    \zero &  -\pdarclen \ba_2^T &   \pdarclen \ba_1^T &               \zero
  \end{bmatrix}
  , &
  \bUbre_1(\ba) &= \frac12
  \begin{bmatrix}
    2 \ba_1^T &    \zero &    \zero &    \zero \\
    2 \ba_2^T &    \zero &    \zero &    \zero \\
    2 \ba_3^T &    \zero &    \zero &    \zero \\
    {}  \zero &    \zero &  \ba_3^T & -\ba_2^T \\
    {}  \zero & -\ba_3^T &    \zero &  \ba_1^T \\
    {}  \zero &  \ba_2^T & -\ba_1^T &    \zero
  \end{bmatrix}
  ,
\end{align}
and
\begin{equation}
  \Doper =
  \begin{bmatrix}
    \bI \pdarclen &         \zero &         \zero &         \zero \\
    {}      \zero & \bI \pdarclen &         \zero &         \zero \\
    {}      \zero &         \zero & \bI \pdarclen &         \zero \\
    {}      \zero &         \zero &         \zero & \bI \pdarclen
  \end{bmatrix}
  .
\end{equation}
As already discussed for assembling the matrix $\KKT$,
we are also required to compute
$\Uoper_2(\stress) = \pd{\coor} (\Boper(\coor) \stress)$,
which is given explicitly by
\begin{equation}
  \Uoper_2(\stress) =
  \Doper^T_{\bgamma} \bUbar_2(\stress) \Doper_{\bgamma} +
  \Doper^T_{\bomega} \bUbre_2(\stress) \Doper_{\bomega},
\end{equation}
where
\begin{align}
  \bUbar_2(\stress) &=
  \begin{bmatrix}
    {}\zero & s_1 \bI & s_2 \bI & s_3 \bI \\
    s_1 \bI &   \zero &   \zero &   \zero \\
    s_2 \bI &   \zero &   \zero &   \zero \\
    s_3 \bI &   \zero &   \zero &   \zero
  \end{bmatrix}
  , &
  \Doper_{\bgamma} &=
  \begin{bmatrix}
    \bI\pdarclen & \zero & \zero & \zero \\
    {}     \zero &   \bI & \zero & \zero \\
    {}     \zero & \zero &   \bI & \zero \\
    {}     \zero & \zero & \zero & \bI
  \end{bmatrix}
\end{align}
as well as
\begin{equation}
  \bUbre_2(\stress) =
  \frac12
  \begin{bmatrix}
    \zero&   \zero&   \zero&   \zero&\zero&   \zero&   \zero&   \zero \\
    \zero&   \zero&   \zero&   \zero&\zero&-s_6 \bI&   \zero& s_5 \bI \\
    \zero&   \zero&   \zero& s_6 \bI&\zero&   \zero&-s_5 \bI&   \zero \\
    \zero&   \zero& s_6 \bI&   \zero&\zero&   \zero&   \zero&-s_4 \bI \\
    \zero&   \zero&   \zero&   \zero&\zero&   \zero&   \zero&   \zero \\
    \zero&-s_6 \bI&   \zero&   \zero&\zero&   \zero& s_4 \bI&   \zero \\
    \zero&   \zero&-s_5 \bI&   \zero&\zero& s_4 \bI&   \zero&   \zero \\
    \zero& s_5 \bI&   \zero&-s_4 \bI&\zero&   \zero&   \zero&   \zero
  \end{bmatrix}
\end{equation}
and
\begin{equation}
  \Doper_{\bomega} =
  \begin{bmatrix}
    {}        \bI &         \zero &         \zero &         \zero \\
    {}      \zero &           \bI &         \zero &         \zero \\
    {}      \zero &         \zero &           \bI &         \zero \\
    {}      \zero &         \zero &         \zero &           \bI \\
    \bI \pdarclen &         \zero &         \zero &         \zero \\
    {}      \zero & \bI \pdarclen &         \zero &         \zero \\
    {}      \zero &         \zero & \bI \pdarclen &         \zero \\
    {}      \zero &         \zero &         \zero & \bI \pdarclen
  \end{bmatrix}
  .
\end{equation}

To perform the spatial discretization of the geometrically exact beam
into two-node finite elements,
we approximate the kinematic fields as well as their
admissible variations with first-order Lagrangian functions.
The adopted numerical scheme for the integration of elemental contributions
is the standard Gauss-Legendre quadrature rule.
As usual, the integrals involving internal terms are computed by means
of a one-point integration scheme that avoids shear locking issues.
Therefore, the evaluation of the kinematic fields at the single Gauss point
is in fact an average of the nodal values,
and their derivatives with respect to the arc length turn out to be
the simplest directed difference of the nodal values. Additionally, we have that $\Boper\leadsto\bB$, $\Uoper_1\leadsto\bU_1$ and $\Uoper_2\leadsto\bU_2$.
Moreover, even for coarse discretizations,
no additional residual stress corrections are necessary.
For an extensive treatment of geometrically exact beams
in a non-\dd~finite element setting, we refer to
\cite{Romero2002, Romero2004, Betsch2002, Gebhardt2019a, Gebhardt2019b}.

Finally, following \cite{Betsch2002}, the mutual orthonormality condition among the directors
is enforced at the nodal level by means of the internal constraint
\begin{equation}
  \cons =
  \frac12 \begin{pmatrix}
    \coori \cdot \coori - 1 \\
    \coorj \cdot \coorj - 1 \\
    \coork \cdot \coork - 1 \\
    2 \coorj \cdot \coork \\
    2 \coori \cdot \coork \\
    2 \coori \cdot \coorj
  \end{pmatrix}
  .
\end{equation}
The associated Jacobian matrix is
\begin{equation}
  \bHcoor =
  \begin{bmatrix}
    \zero & \coori^T & \zero    & \zero \\
    \zero & \zero    & \coorj^T & \zero \\
    \zero & \zero    & \zero    & \coork^T \\
    \zero & \zero    & \coork^T & \coorj^T \\
    \zero & \coork^T & \zero    & \coori^T \\
    \zero & \coorj^T & \coori^T & \zero
  \end{bmatrix}
\end{equation}
and $\bV(\bnu) = \pd{\coor} (\bHcoorT \bnu)$ is simply
\begin{equation}
  \bV(\bnu) =
  \begin{bmatrix}
    \zero &     \zero &     \zero &     \zero \\
    \zero & \nu_1 \bI & \nu_6 \bI & \nu_5 \bI \\
    \zero & \nu_6 \bI & \nu_2 \bI & \nu_4 \bI \\
    \zero & \nu_5 \bI & \nu_4 \bI & \nu_3 \bI
  \end{bmatrix}
  ,
\end{equation}
which can be interpreted as an additional stiffness
due to the presence of the internal constraint.

The null-space basis corresponding to the internal constraint
at the nodal level can be built by visual inspection of the Jacobian matrix as
\begin{equation}
  \nullproj =
  \begin{bmatrix}
    \bI & \zero & \zero & \zero \\
    \zero & \coorihat & \coorjhat & \coorkhat
  \end{bmatrix}^T
  ,
\end{equation}
where the algebraic operator $\widehat{(\rdot)}$ emulates the cross product.

The matrix $\bW_1(\ba) = \pd{\coor} (\nullproj \ba)$ is
\begin{equation}
  \bW_1(\ba) =
  -\begin{bmatrix}
    \zero & \zero & \zero & \zero \\
    \zero & \baihat & \bajhat & \bakhat
  \end{bmatrix}
  ,
\end{equation}
and the matrix $\bW_2(\bb) = \pd{\coor} (\nullprojT \bb)$ is
\begin{equation}
  \bW_2(\bb) =
  \begin{bmatrix}
    \zero &     \zero &     \zero & \zero \\
    \zero & \bbtwohat &     \zero & \zero \\
    \zero &     \zero & \bbtwohat & \zero \\
    \zero &     \zero &     \zero & \bbtwohat
  \end{bmatrix}
  ,
\end{equation}
where $\bb = (\bb_1^T, \bb_2^T) \in \reals^6$
with $\bb_1$ and $\bb_2$ in $\reals^3$.

The formulation of constraints that can be related to the usual types of boundary conditions (rigid support, simple support, movable support \ia),
can be achieved by linear equations in the nodal variables and is thus straightforward. In the scope of this work, we therefore omit their systematic presentation and refer for further details to \cite{Hente2019,Gebhardt2019a}.
With all vectors and matrices derived in this subsection at hand, the construction of the equations corresponding to the optimality conditions and their derivatives is straightforward.

%% file: 5examples.tex
\section{Numerical examples}

\begin{figure}[tp]
  \centering
  \includegraphics[width=0.465\textwidth]{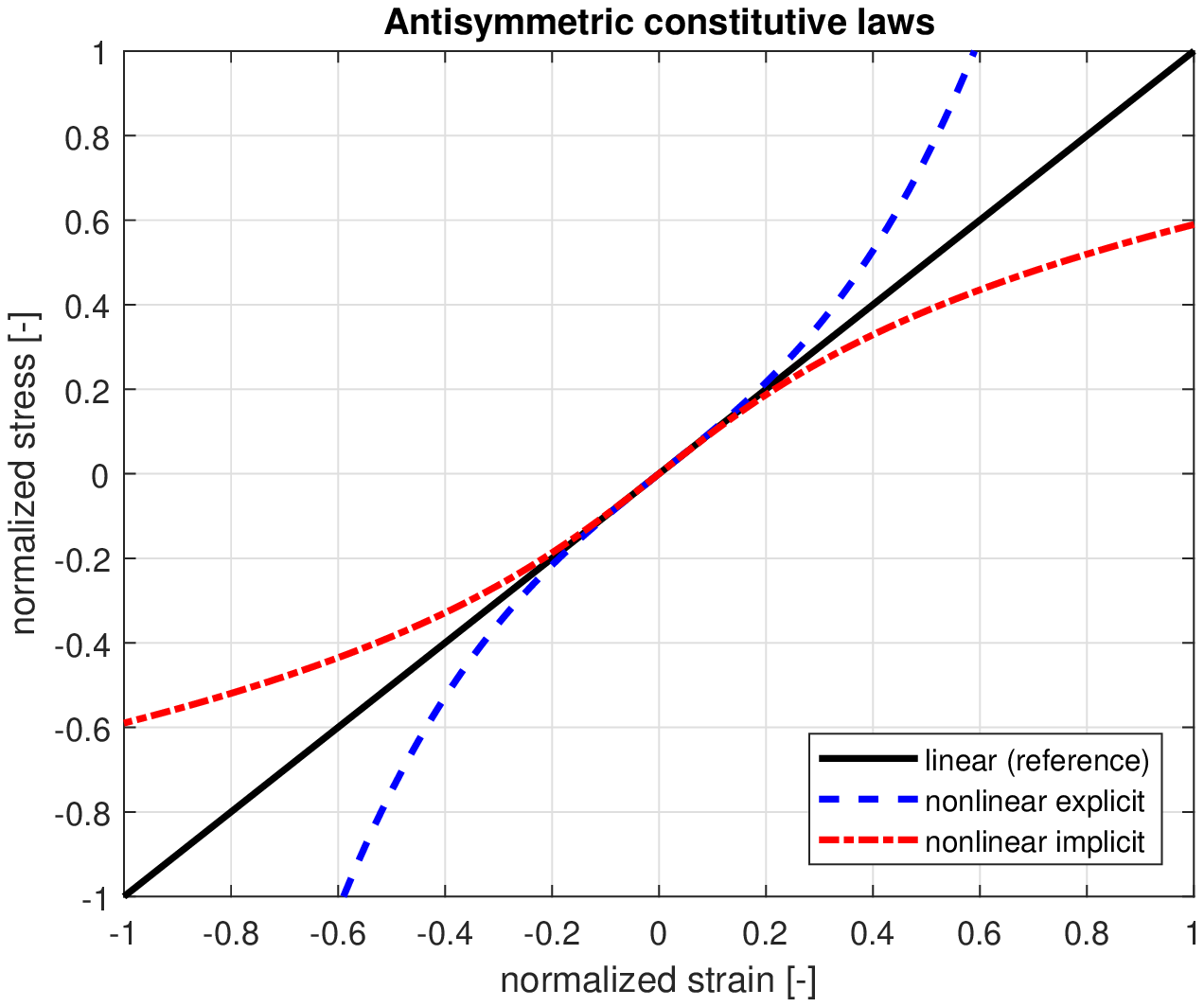}\hfill
  \includegraphics[width=0.465\textwidth]{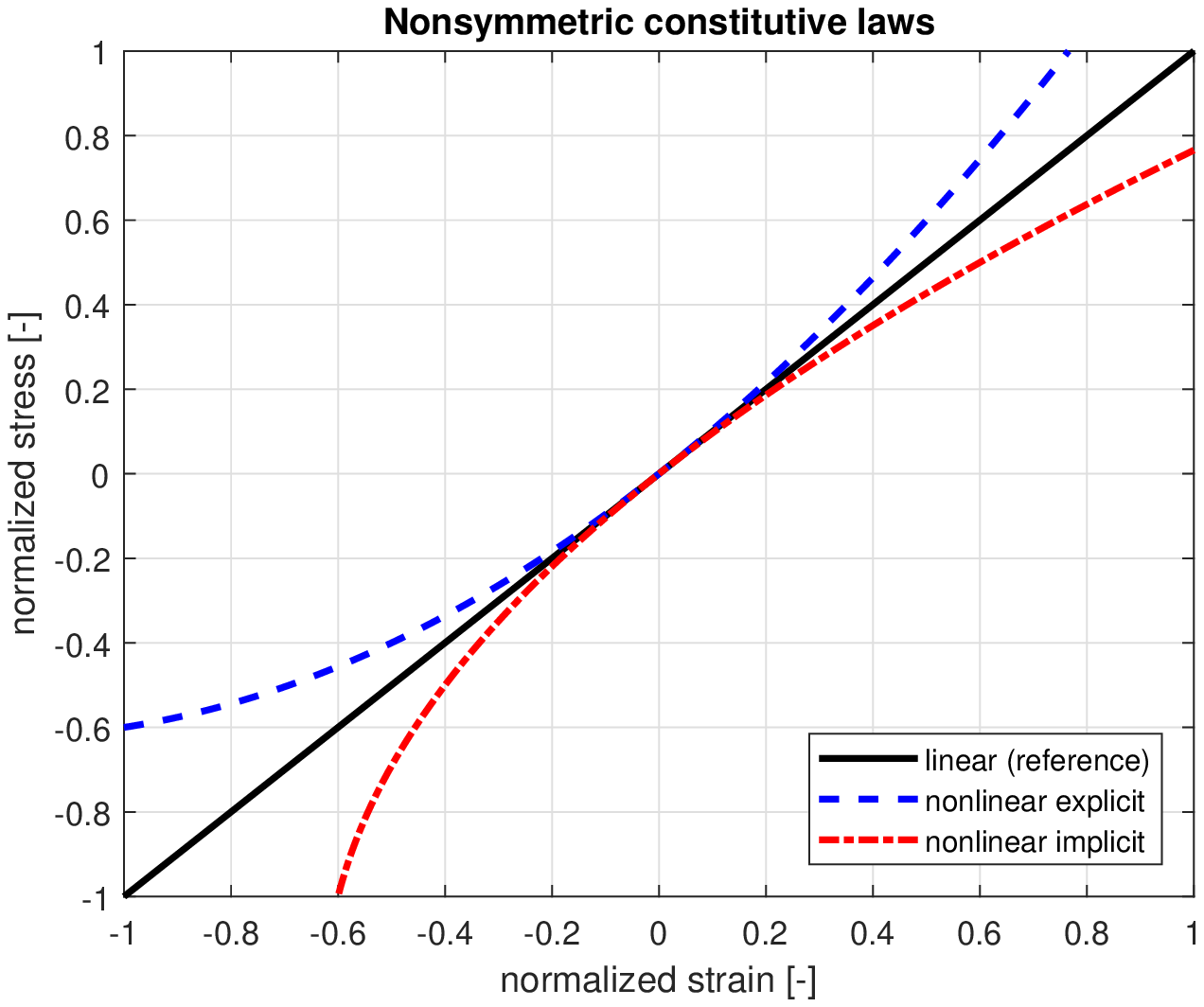}%
  \caption{Nonlinear constitutive laws
    with explicitly or implicitly defined stress.}
  \label{fig:materials}
\end{figure}

In this section, we present three numerical examples of increasing difficulty to show the potential of the proposed approximate NLP for \DDCM. Specifically, we consider its specialization to the geometrically exact beam model. The first example presents a verification of the proposed formulation, taking the underlying three-director based standard FE formulation for geometrically exact beams as a reference (see \cite{Romero2002, Betsch2002, Gebhardt2019a, Gebhardt2019b} for further details). The reference formulation is equipped with the simplest linear constitutive law. In the second example, we investigate the numerical behavior of our framework for antisymmetric and nonsymmetric forms of nonlinear constitutive laws with explicitly defined stress resultants, \ie $\stress=\stress(\strain)$. Their basic nonlinear behavior in the phase space is illustrated in Figure \ref{fig:materials}. In the third example, we investigate the numerical behavior of our framework for nonlinear constitutive laws in antisymmetric and nonsymmetric forms, in which the stress resultants are implicitly defined, \ie $\stress\neq\stress(\strain)$, see again Figure \ref{fig:materials}. We note that the antisymmetric case is motivated by experimental results with state diagrams that render ``S''-shaped graphs. For representing such mechanical behavior, we find the formulation of a general constitutive manifold more natural.

It is clear that problem types represented by the first and second examples can also be addressed by any standard FE formulation. The problem type represented by the third example can be naturally addressed in the context of our approximate NLP approach, but not by standard FEM, and therefore illustrates a possible advantage of our approach over standard FEM. We believe that this advantage of our approach opens a new dimension towards the analysis of materials represented by non-convex constitutive laws such as metamaterials, multiscale materials, \etc. Whether the materials considered are physically feasible or not is beyond the scope of the current work and deserves an extensive investigation regarding their functional structure and identification \ia. In this context, techniques like manifold learning and nonlinear dimensionality reduction seem to be very promising approaches.

% Dominik: Der letzte Satz hängt bisschen in der Luft. Vielleicht wäre es gut, ein zwei Sätze hinzuzufügen, die die Verbindung zu unserer Methode besser beschreibt.

\begin{figure}[tp]
  \centering
  \includegraphics[height=0.35\textheight]{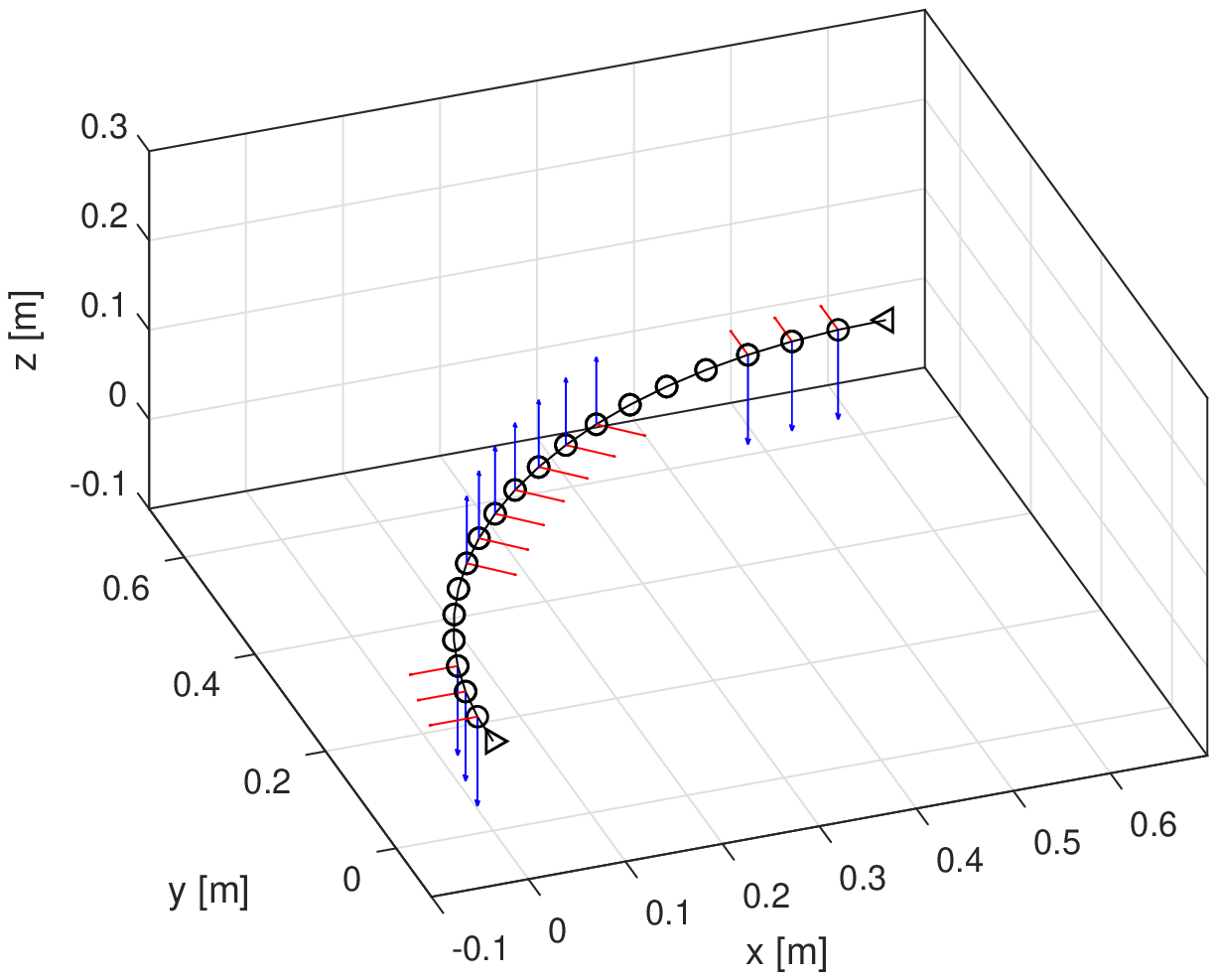}
  \caption{Finite element representation of the beam structure.
    Circles indicate nodes with internal constraints only. Triangles indicate nodes that are fully fixed.
    Blue and red arrows denote forces along the vertical direction $z$ and in the horizontal $x$-$y$-plane, respectively.}
  \label{fig:resBeamStructure}
\end{figure}

All the three examples are built on a curved beam structure whose geometry is described by a quarter of a circular arc with a total arc length of $1\,\textrm{m}$, which corresponds to a radius of $\frac{2}{\pi}\,\textrm{m}$. Both ends are fully fixed and the structure is uniformly discretized into $21$ finite elements. Figure \ref{fig:resBeamStructure} shows the finite element representation, the boundary conditions and loads applied. The first node is located at the position $(0,0,0)\,\textrm{m}$. The nodes $2$--$4$ are loaded with vertical nodal forces $(0,0,-20)\,\textrm{N}$ and horizontal nodal forces $(-10,0,0)\,\textrm{N}$. For the nodes $8$--$14$ we have $(0,0,15)\,\textrm{N}$ and $(7.5,-7.5,0)\,\textrm{N}$. Lastly, for the nodes $18$--$20$ we have $(0,0,-20)\,\textrm{N}$ and $(0,10,0)\,\textrm{N}$. This setting has been chosen, because the action of combined spatial loads creates complex strain-stress states and the geometrically nonlinear response is apparent at first glance. Moreover, for all numerical examples the load was applied in one single load step and a relative error-based tolerance of $10^{-12}$ has been set for the Newton iteration.

\subsection{Verification}

% In this first numerical example, we present a verification for the approximate NLP approach. Clearly, it does not verify the approach thoroughly, but it is still good enough for the purpose of our work.
% Dominik: Why do you weaken your line of arguments before you even started?

As indicated above, we first consider the beam structure described above and the simplest linear material law, defined as
\begin{equation}
(\bg^\sharp(\strainrec,\stressrec))^i =\check{s}^i-a^{ii}\check{e}_i=0
\qquad\textrm{or alternatively as}\qquad
(\bg_\flat(\strainrec,\stressrec))_i =\check{e}_i-a_{ii}\check{s}^i=0,
\end{equation}
\ifcase0
with $a^{ii}=a^{-1}_{ii}$ and the (nonphysical) values $a^{11}=a^{22}= 75\,\mathrm{N}$, $a^{33}= 100\,\mathrm{N}$, $a^{44}=a^{55}= 100\,\mathrm{Nm^2}$ and $a^{66} = 200\,\mathrm{Nm^2}$. The weight matrix $\weight$ is defined as the identity.
\or
with $a^{ii}=a^{-1}_{ii}$ and the (nonphysical) values $a^{11}=a^{22}= 7.5 \times 10\,\mathrm{N}$, $a^{33}= 1.0 \times 10^2\,\mathrm{N}$, $a^{44}=a^{55}= 1.0 \times 10^2\,\mathrm{Nm^2}$ and $a^{66} = 2.0 \times 10^2\,\mathrm{Nm^2}$. The weight matrix $\weight$ is defined as the identity.
\fi

\begin{figure}[tp]
  \centering
  \includegraphics[height=0.33\textheight]{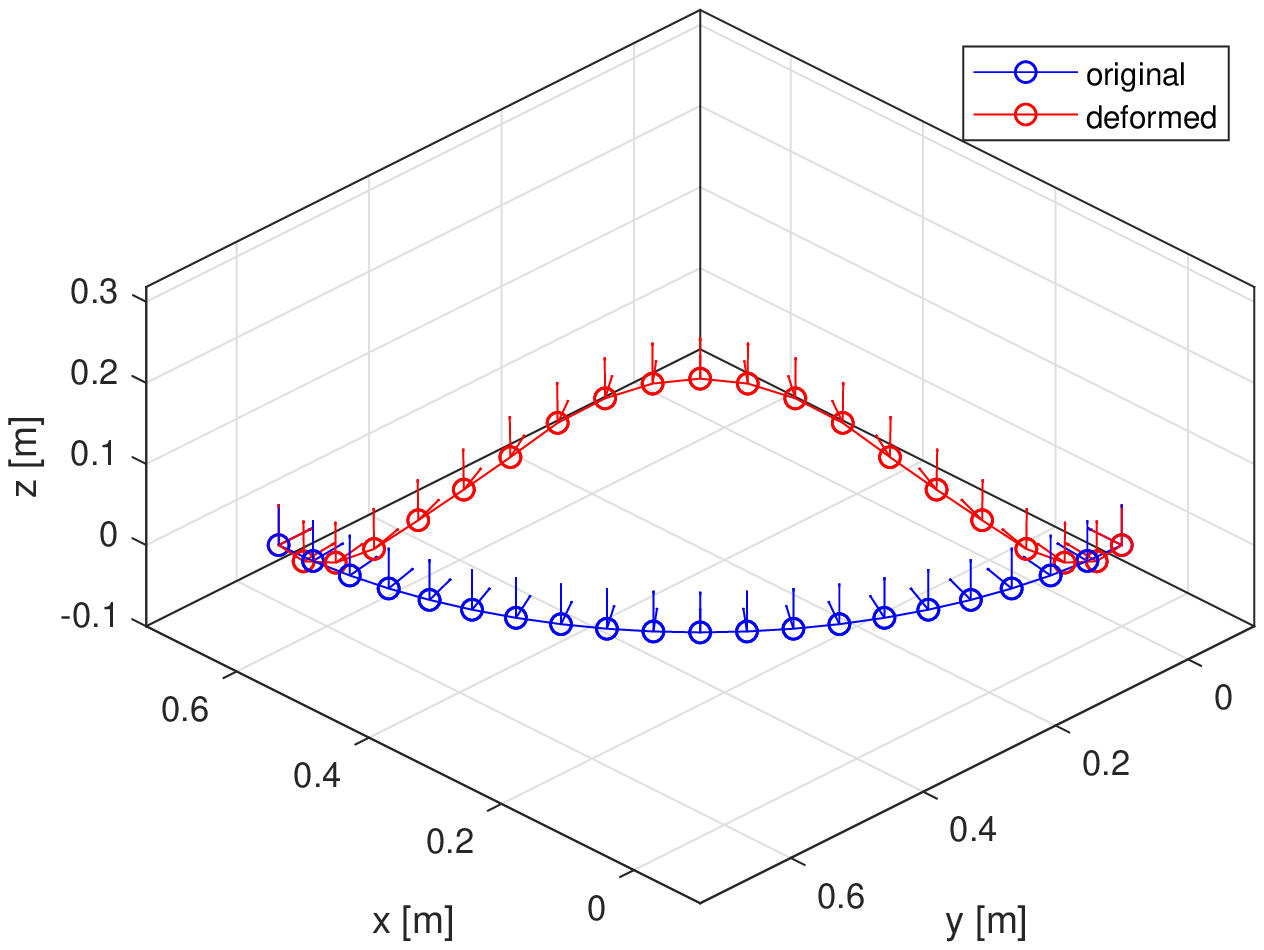}
  \caption{Verification (first example): original and deformed configurations.}
  \label{fig:verorivsdef}
\end{figure}

Simulations were carried out with standard FEM as well as with our approximate NLP approach. The Newton iteration was started with the stress-free state and zero Lagrange multipliers in both cases. Figure \ref{fig:verorivsdef} shows the original and deformed unscaled configurations computed with the approximate NLP approach, illustrating the nonlinear response.
Both approaches require five iterations to find a solution. After the first two iterations, we observed quadratic converge with both methods, which we think is a strong indicator of quality regarding the correctness of the implementations.
% \MCS[In SFEF, a nonlinear system of equations is solved: which?] Later!

\ifcase0
\begin{table}[tp]
  \centering
  \def\-{$-$}
  \def\+{\hphantom{$-$}}
  \begin{tabular}{ccccc}
    \toprule
    comp. & $\coorz$ & $\coori$ & $\coorj$ & $\coork$ \\
    \midrule
    $x$ & 0.30620366 & \+0.00761486 & \+0.70706900 & 0.70710999 \\
    $y$ & 0.33041646 & \-0.00761485 & \-0.70706899 & 0.70711000 \\
    $z$ & 0.21515427 & \+0.99994201 & \-0.01076908 & 0.00000000 \\
    \midrule
    $x$ & 0.30620367 & \+0.00761487 & \+0.70706900 & 0.70711000 \\
    $y$ & 0.33041647 & \-0.00761485 & \-0.70706899 & 0.70711000 \\
    $z$ & 0.21515428 & \+0.99994201 & \-0.01076909 & 0.00000000 \\
    \bottomrule
  \end{tabular}
  \caption{Nodal variables of the $11^\textrm{th}$ node at the deformed
    configuration; standard FEM (top) vs.\ approximate NLP (bottom).}
  \label{tab:sfn11}
\end{table}
\or
\begin{table}[tp]
  \centering
  \def\-{$-$}
  \def\+{\hphantom{$-$}}
  \begin{tabular}{ccccc}
    \toprule
    comp. & $\coorz$ & $\coori$ & $\coorj$ & $\coork$ \\
    \midrule
    $x$ & 0.30620366 & \+0.00761486 & \+0.70706900 & 0.70710999 \\
    $y$ & 0.33041646 & \-0.00761485 & \-0.70706899 & 0.70711000 \\
    $z$ & 0.21515427 & \+0.99994201 & \-0.01076908 & 0.00000000 \\
    \bottomrule
  \end{tabular}
  \caption{Standard FEM: Nodal variables of the
    $11^\textrm{th}$ node at the deformed configuration.}
  \label{tab:sfn11}
\end{table}

\begin{table}[tp]
  \centering
  \def\-{$-$}
  \def\+{\hphantom{$-$}}
  \begin{tabular}{ccccc}
    \toprule
    comp. & $\coorz$ & $\coori$ & $\coorj$ & $\coork$ \\
    \midrule
    $x$ & 0.30620367 & \+0.00761487 & \+0.70706900 & 0.70711000 \\
    $y$ & 0.33041647 & \-0.00761485 & \-0.70706899 & 0.70711000 \\
    $z$ & 0.21515428 & \+0.99994201 & \-0.01076909 & 0.00000000 \\
    \bottomrule
  \end{tabular}
  \caption{Approximate NLP: Nodal variables of the
    $11^\textrm{th}$ node at the deformed configuration.}
  \label{tab:aolnlpn11}
\end{table}
\fi

\begin{table}[tp]
  \centering
  \def\-{$-$}
  \def\+{\hphantom{$-$}}
  \begin{tabular}{ccccc}
    \toprule
    comp. &
    $\strain$ (std.~FEM) & $\strain$ (approx.~NLP) &
    $\stress$ (std.~FEM) & $\stress$ (approx.~NLP) \\
    \midrule
    $1$ & \+0.50365604 & \+0.50365604 & \+37.77420296 & \+37.77420295 \\
    $2$ & \+0.34693364 & \+0.34693365 & \+26.02002346 & \+26.02002347 \\
    $3$ & \-0.02990629 & \-0.02990630 &  \-2.99062981 &  \-2.99062982 \\
    $4$ & \+0.04172209 & \+0.04172209 &  \+4.17220924 &  \+4.17220925 \\
    $5$ & \-0.05823386 & \-0.05823387 &  \-5.82338675 &  \-5.82338675 \\
    $6$ & \+0.00749371 & \+0.00749371 &  \+1.49874166 &  \+1.49874166 \\
    \bottomrule
  \end{tabular}
  \caption{Stress-strain states
    of the $8^\textrm{th}$ element at the deformed configuration.}
  \label{tab:sfe8andaolnlpe8}
\end{table}

\iffalse
\begin{table}[h!]
	\centering{}
	\begin{tabular}{|c|c|c|}
		\hline
		comp. & $\strain$ & $\stress$\\
		\hline
		$1$ & 0.50365604 & 37.77420295 \\
		$2$ & 0.34693365 & 26.02002347 \\
		$3$ &-0.02990630 &-2.99062982 \\
		$4$ & 0.30620367 & 4.17220925 \\
		$5$ &-0.05823387 &-5.82338675 \\
		$6$ & 0.00749371 & 1.49874166 \\
		\hline
	\end{tabular}
	\caption{Approximate NLP: stress-strain states of the $8^\textrm{th}$ element at the deformed configuration.}
	\label{tab:aolnlpe8}
\end{table}
\fi

\ifcase0
Table \ref{tab:sfn11} presents the position vector and directors corresponding to the $11^\textrm{th}$ node, the central one, at the deformed configuration computed with the standard FEM (top) and the same information computed with the approximate NLP (bottom). As it can be observed, the agreement among the results is excellent. Table \ref{tab:sfe8andaolnlpe8} presents the stress resultants corresponding to the $8^\textrm{th}$ element at the deformed configuration computed with the standard FEM and with the approximate NLP. Once again, the agreement is excellent.
\or
Table \ref{tab:sfn11} presents the position vector and directors corresponding to the $11^\textrm{th}$ node, the central one, at the deformed configuration computed with the standard FEM. Table \ref{tab:aolnlpn11} presents the same information computed with the approximate NLP. As it can be observed, the agreement among the results is excellent. Table \ref{tab:sfe8andaolnlpe8} presents the stress resultants corresponding to the $8^\textrm{th}$ element at the deformed configuration computed with the standard FEM and with the approximate NLP. Once again, the agreement is excellent.
\fi

\newcommand\sixpack[3]{%
  \begin{figure}[tp]
    \centering
    \includegraphics[width=0.465\textwidth]{#1_1}\hfill
    \includegraphics[width=0.465\textwidth]{#1_4}\\[2.45ex]
    \includegraphics[width=0.465\textwidth]{#1_2}\hfill
    \includegraphics[width=0.465\textwidth]{#1_5}\\[2.45ex]
    \includegraphics[width=0.465\textwidth]{#1_3}\hfill
    \includegraphics[width=0.465\textwidth]{#1_6}%
    \caption{#2: Distribution of the stress resultant components
      1, 2, 3 (left) and 4, 5, 6 (right).
      The components 1, 2, 3 correspond to the cross sectional
      force vector per unit length;
      the components 4, 5, 6 correspond to the cross sectional
      moment vector per unit length.}
    \label{fig:#3}
  \end{figure}
}
\sixpack{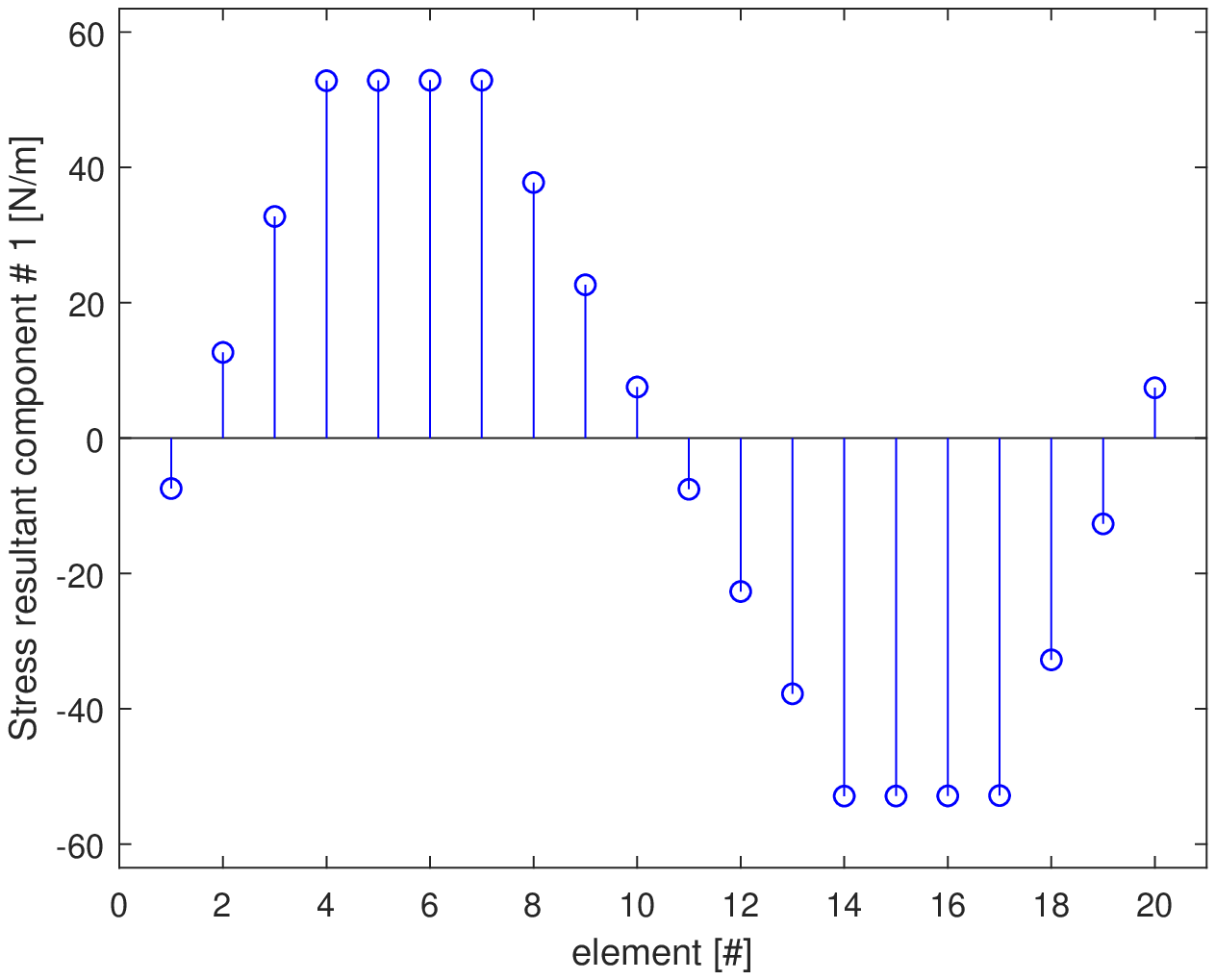}{Verification (first example)}{verstress}
Figure~\ref{fig:verstress} shows the distribution of all stress resultants (with the sign of the local frame) over the elements that were consecutively arranged along the arc length of the beam. These results were computed with the approximate NLP. We found again an excellent agreement with the results computed with the standard FEM. For the sake of brevity, the results obtained with the standard FEM are omitted here. We note that due to the symmetry of the structure, boundary conditions and applied loads, the resultant transversal forces are antisymmetric and the resultant transversal moments are symmetric. The resultant axial force is symmetric and the resultant torsion moment is antisymmetric.

\subsection{Explicit stress definition}

In this second numerical example, we consider nonlinear constitutive laws of the following form
\begin{equation}
(\bg^\sharp(\strainrec,\stressrec))^i =\check{s}^i-a^{ii}\check{e}_i-\frac{b^{ii}}{2}(\check{e}_i)^2-\frac{c^{iiii}}{3}(\check{e}_i)^3,
\end{equation}
where the stress resultants are defined \textit{explicitly} in terms of the strain measures.
Notice that more elaborate equations of non-polynomial form are possible as well, and will be subject of future investigations. For the coefficients $a^{ii}$, we use the same (nonphysical) values as in the verification example. Moreover, we set $b^{iii} = 0.85 a^{ii}$ and $c^{iiii} = a^{ii}$.

\begin{figure}[tp]
  \centering
  \includegraphics[height=0.33\textheight]{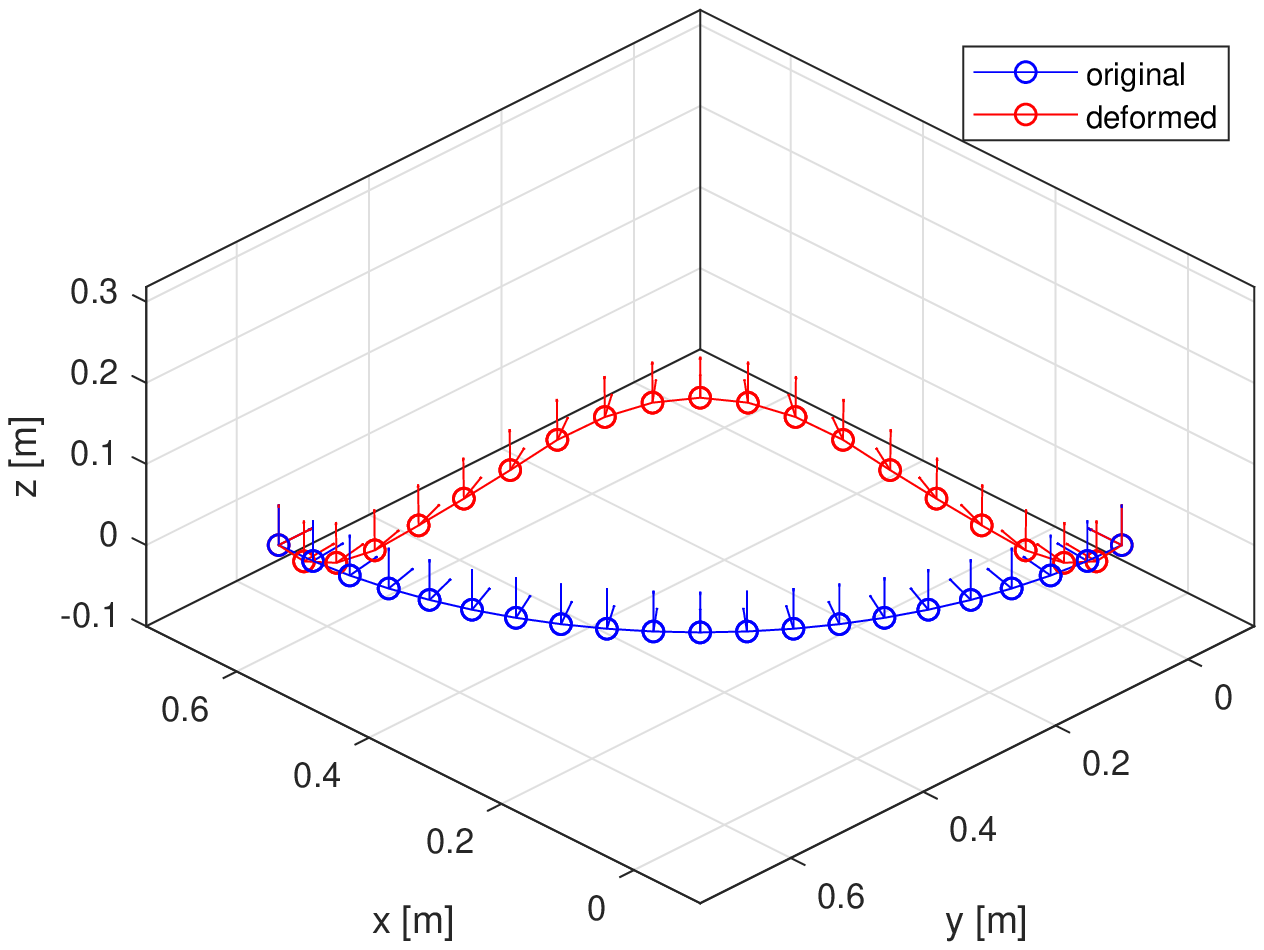}
  \caption{Antisymmetric explicit stress definition:
    Original and deformed configurations.}
  \label{fig:ex2c1}
\end{figure}

Next, we consider two cases, where all computations are performed with the approximate NLP. In a first case, we consider only linear and cubic contributions. Such a setting preserves the antisymmetric feature of the material response. Figure \ref{fig:ex2c1} shows the original and deformed configurations. Due to the symmetry of the beam geometry, the boundary conditions and the applied loads in combination with the antisymmetry of the material response, the resulting deformed configuration is perfectly symmetric. The global minimizer was found in six iterations. As observed in the verification case, the distributions of the resultant transversal forces are antisymmetric and the distributions of the resultant transversal moments are symmetric. The distribution of the resultant axial force is symmetric and the distribution of the resultant torsion moment is antisymmetric. Both exhibit the expected mechanical behavior and are not shown here.

\begin{figure}[tp]
  \centering
  \includegraphics[height=0.33\textheight]{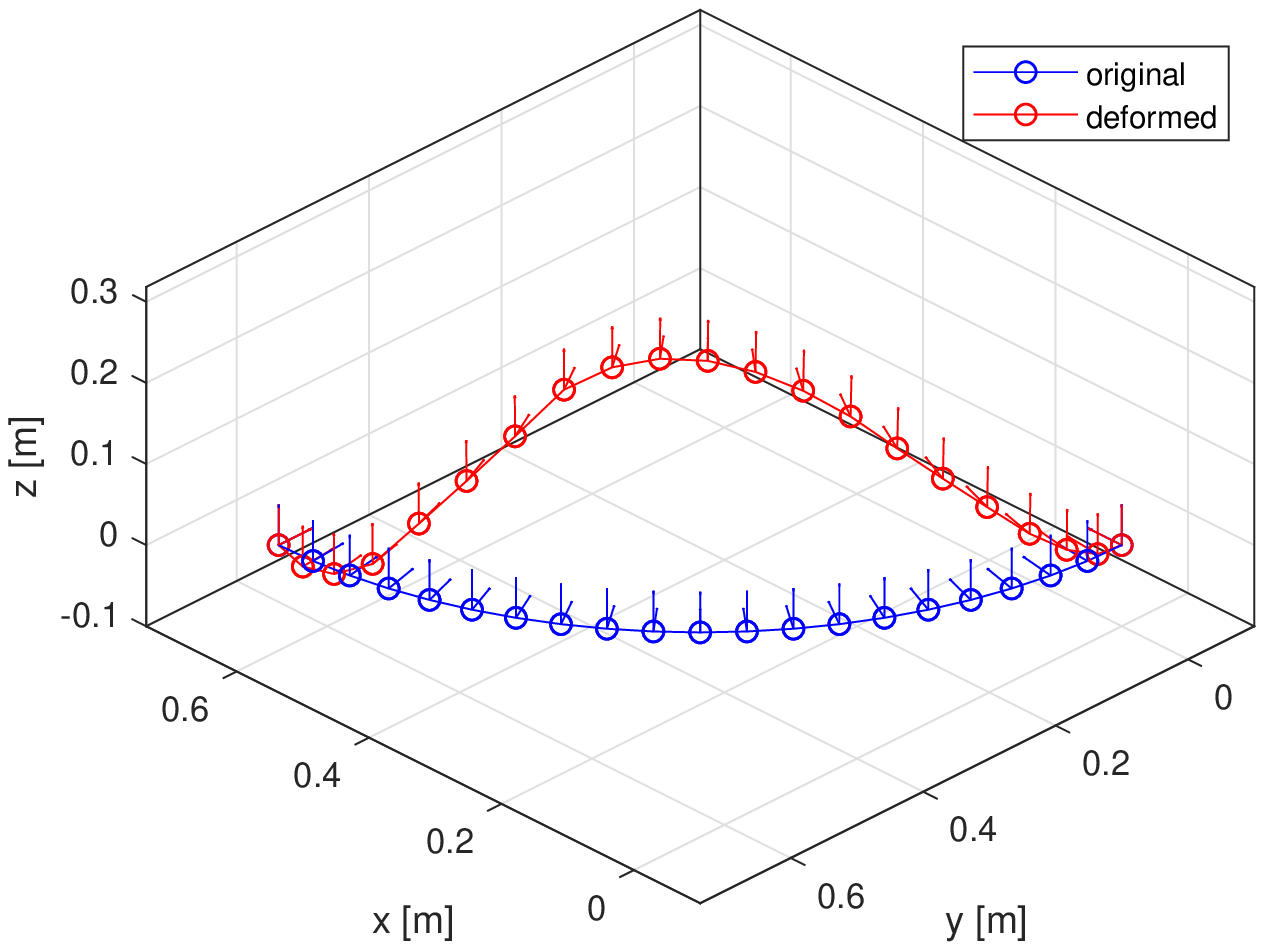}
  \caption{Nonsymmetric explicit stress definition:
    Original and deformed configurations.}
  \label{fig:ex2c2}
\end{figure}

In the second case, we consider only linear and quadratic contributions. That setting destroys the antisymmetric feature of the material response and renders a nonsymmetric response of the structure. Figure \ref{fig:ex2c2} shows the original and deformed configurations. Even though the geometry, boundary conditions and applied loads are symmetric, the resulting deformed configuration is nonsymmetric due to the asymmetry of the material response. The global minimizer was found in seven iterations. In contrast to the previous case, no special symmetry in the distribution of the stress resultants is to be found. Figure~\ref{fig:ex2c2stress} clearly illustrates this feature.

\sixpack{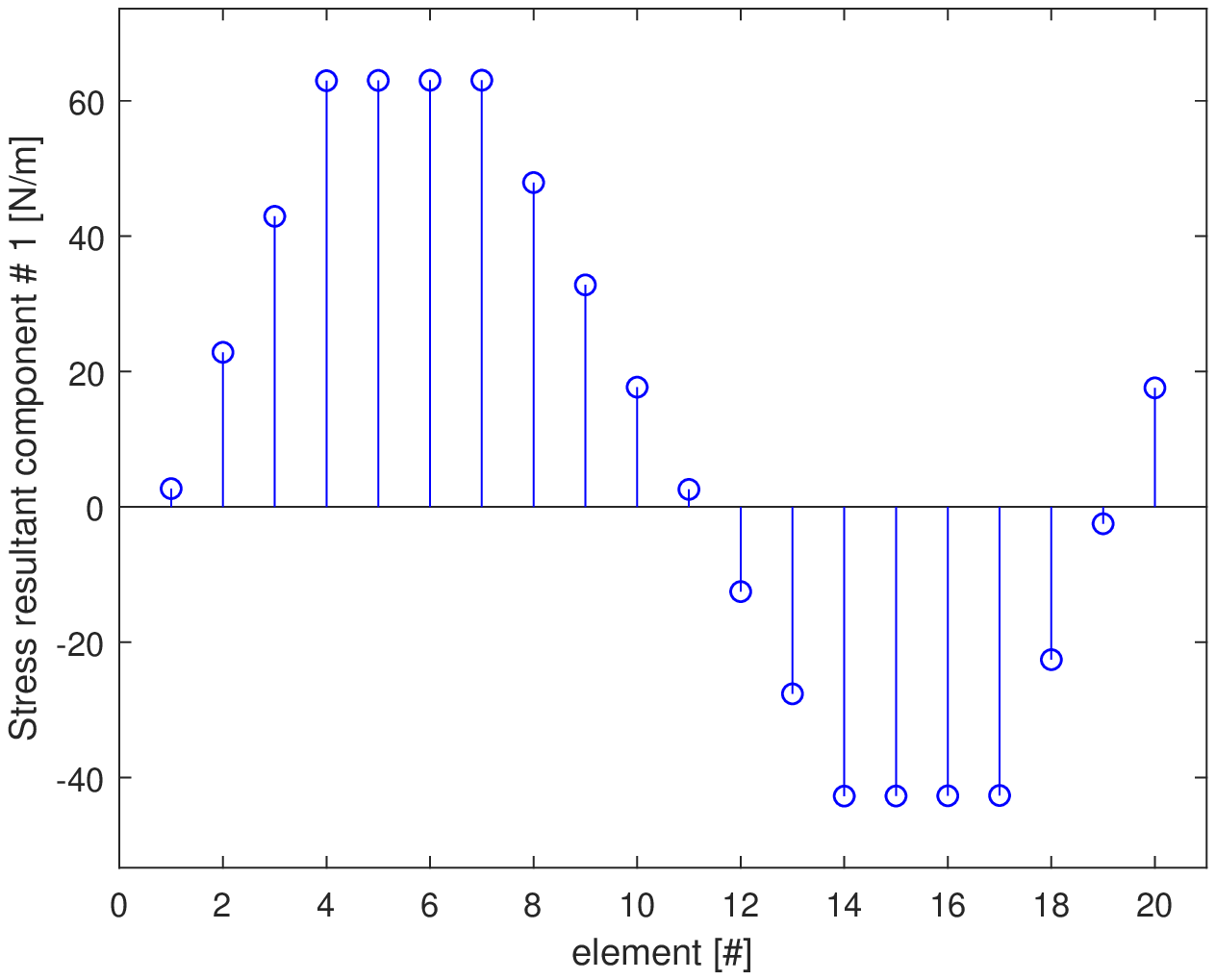}{Nonsymmetric explicit stress definition}{ex2c2stress}

\subsection{Implicit stress definition}

In this third numerical example, we consider nonlinear constitutive laws of the following form
\begin{equation}
(\bg_\flat(\strainrec,\stressrec))_i =\check{e}_i-a_{ii}\check{s}^i-\frac{b_{iii}}{2}(\check{s}^i)^2-\frac{c_{iiii}}{3}(\check{s}^i)^3.
\end{equation}
where the stress resultants are defined \textit{implicitly} in terms of the strain measures.
Once again, more elaborate equations of non-polynomial form are possible as well and will be subject of future investigations. For the coefficients $a_{ii}$, we consider the same (nonphysical) values as in the verification example. Moreover, we set $b_{iii} = 0.015 a_{ii}$ and $c_{iiii} = 0.0005 a_{ii}$. These values are chosen in order to avoid abrupt growth of the stress resultants when the coefficients vary slightly. This is critical especially when the stress resultants take (absolute) values larger than one. We note that in the case of their explicit definition considered in the previous examples, this sensitivity does not occur.

\begin{figure}[tp]
  \centering
  \includegraphics[height=0.33\textheight]{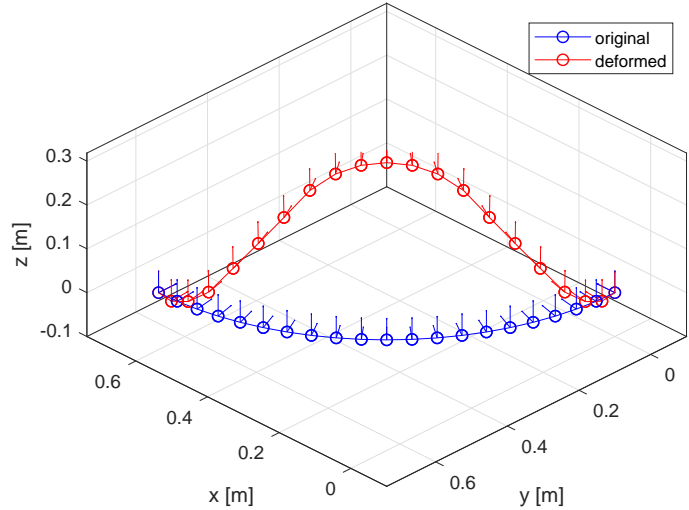}
  \caption{Antisymmetric implicit stress definition:
    Original and deformed configurations.}
  \label{fig:ex3c1}
\end{figure}

\begin{figure}[tp]
  \centering
  \includegraphics[height=0.33\textheight]{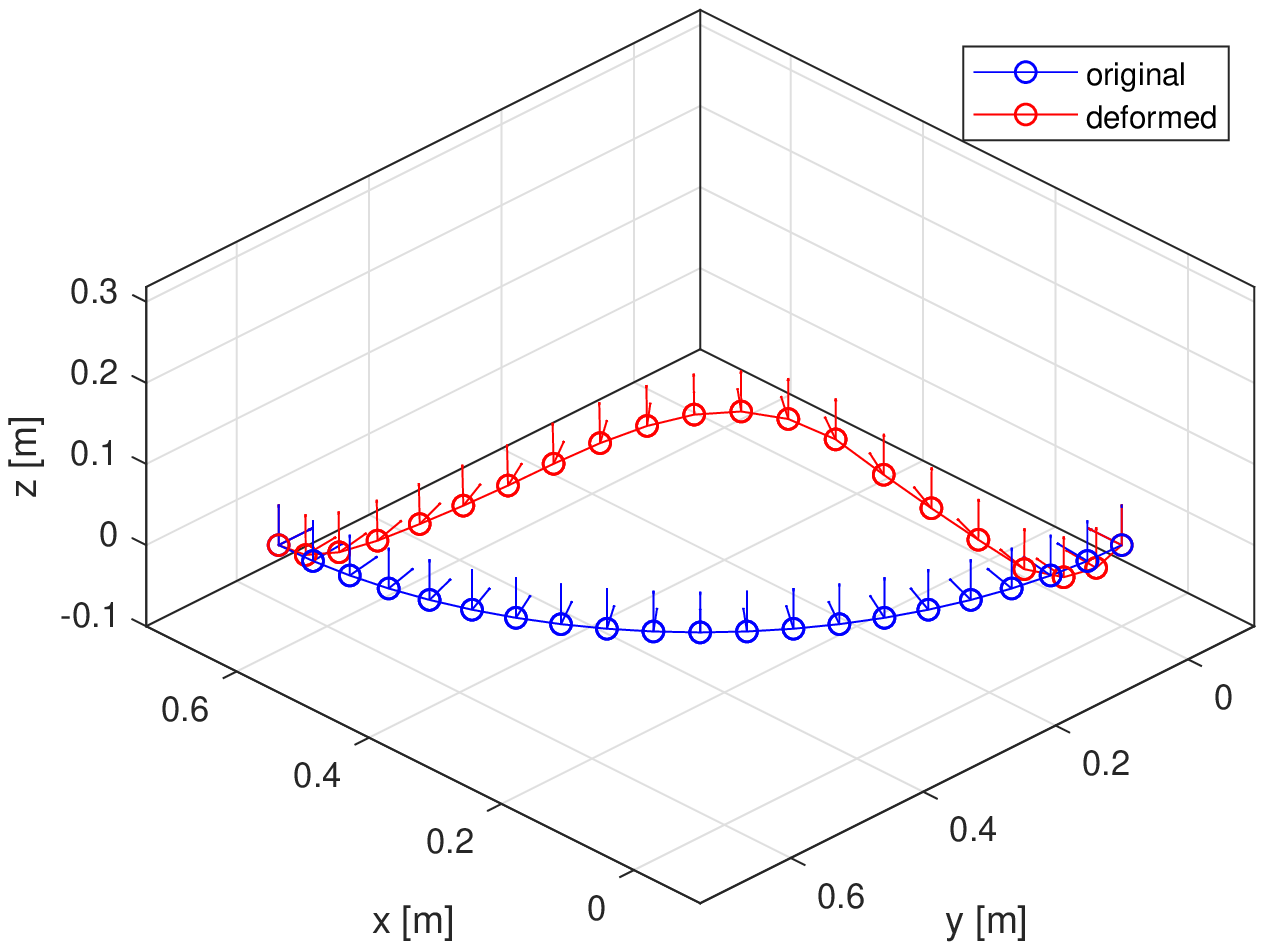}
  \caption{Nonsymmetric implicit stress definition:
    Original and deformed configurations.}
  \label{fig:ex3c2}
\end{figure}

Again we consider two cases. First, we only take into account the linear and cubic contributions such that the antisymmetry of the material response is preserved. Figure \ref{fig:ex3c1} shows the original and deformed configurations. As discussed above, the resulting deformed configuration is perfectly symmetric. The global minimizer was found in five iterations. The distributions of the resultant forces and moments exhibit the expected behavior discussed above. For the sake of brevity, we do not show the complete set of distributions. In the second case, we take into account the linear and quadratic contributions. This setting impedes the antisymmetry of the material response and leads to a nonsymmetric response of the structure. Figure \ref{fig:ex3c2} shows the original and deformed configurations. As discussed above, the resulting deformed configuration is nonsymmetric. The global minimizer was found in six iterations. No special symmetry in the distribution of stress resultants is to be found, see Figure~\ref{fig:ex3c2stress}.

% Dominik: When doing copy-paste of result sections, I would either modify or condense the essential information...

\sixpack{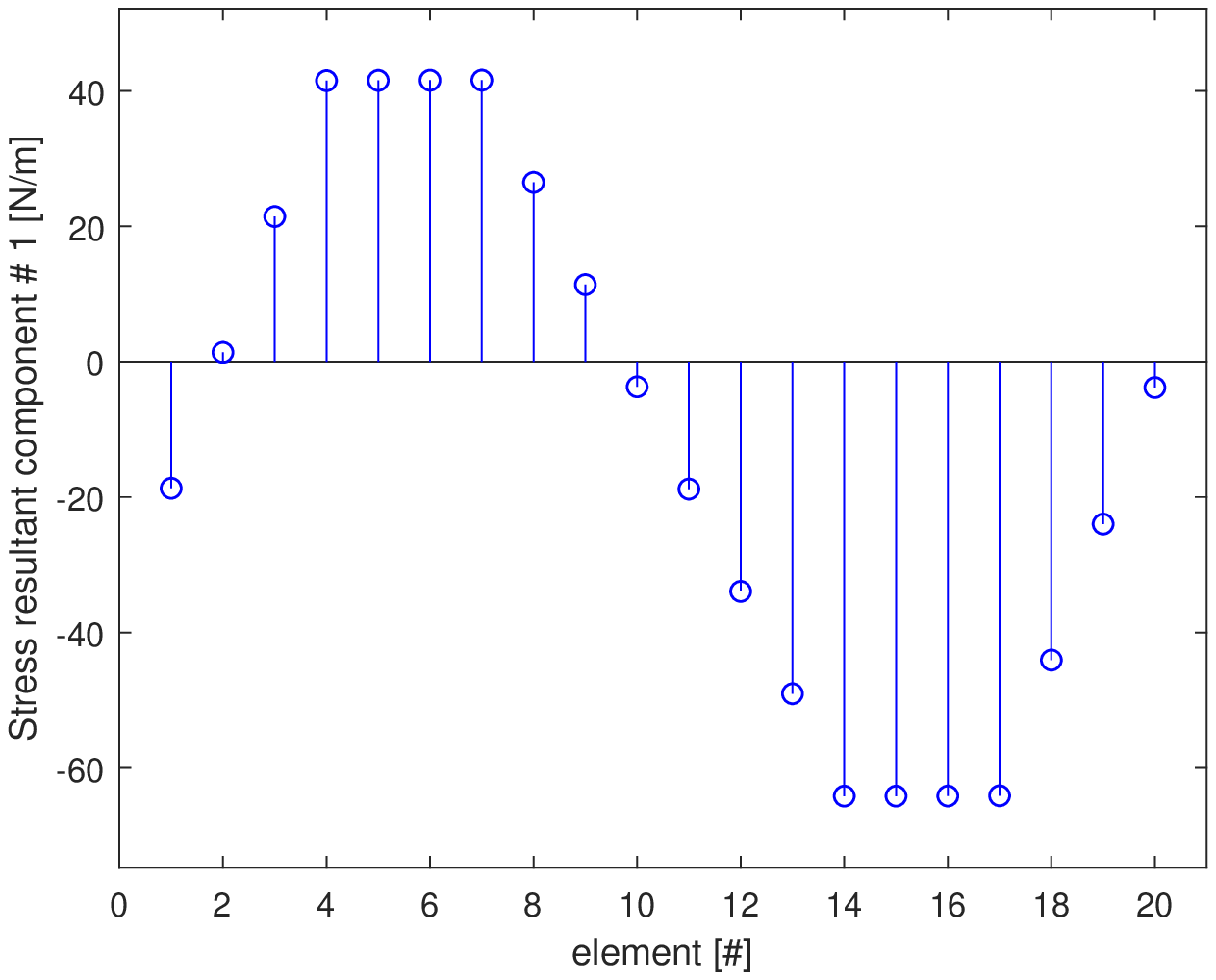}{Nonsymmetric implicit stress definition}{ex3c2stress}

We would like to emphasize again that the results here can only be obtained with our new approximate NLP approach, and therefore a comparison to standard (displacement-based) FEM cannot be provided.

%% file: 6conclusion.tex
\section{Conclusions}

In this work, we presented %for the first time 
an approximate nonlinear optimization problem for \DDCM~that enables us to handle: \textit{i}) kinematic constraints; and, \textit{ii}) materials whose stress-strain relationship can be implicitly approximated. Therefore, our method does not rely on any special functional structure. These features open a new path with respect to the two existing main approaches. These are \textit{i}) discrete-continuous optimization problems that rely on unprocessed data, and \textit{ii}) the identification of energy functions from data by manifold learning techniques. Moreover, our new method combines some strengths of both aforementioned approaches and mitigates some of their weaknesses, such as relatively poor convergence properties and strong sensitivity to the scattering of the data set for the first one, and a special functional structure limited to the explicit definition of stress for the second one. 

The mathematical framework of our approach was completely derived and presented in a self-contained description. In a first step, we applied our new approach to a simple truss element to illustrate the underlying mathematical machinery. In a second step, we successfully applied our method to a geometrically exact beam formulation. We emphasize that the extension of the \DDCM~paradigm to geometrically exact beam elements is by itself new, and to the best of our knowledge, there exists no work in the literature that can deal with this type of elements in the same generality in the context of \DDCM.
%To the best of our knowledge, there is no work in the literature that \reformulate{deals with this particular kind of structural elements} in the context of \DDCM. \MCS[Better ``that can deal with this kind of elements in the same generality'' or similar?] 
The convergence properties for the approximate nonlinear optimization problem are very favorable, as they are comparable to the convergence properties of the underlying standard finite element formulation. This is due to the fact that the approach relies on the Newton method, as opposed to other possible choices such as meta-heuristic methods. In our numerical tests, we observed that our approach is robust, efficient and versatile. 
%\textcolor{red}{DS: I think we should add one or two sentences with key observations from the results described in Section 5 here.}
We therefore think it constitutes a promising starting point for future studies in this direction.
% Dominik: Last sentence -> Did you mean this when you added the work "promising method"?

Further research will address, for instance, the robustness of the approach, which can be improved by using more sophisticated 
% though standard -> Dominik: why do we add a negative connotation? Many engineers will consider it favorable that one can fall back on already established technology that will work for sure...
nonlinear optimization algorithms that can handle inequality constraints and that include a globalization by line-search or trust-region techniques. Advanced techniques for the off-line identification of constitutive manifolds such as manifold learning can be considered.
% and the embedding of the current approach into a multiscale analysis may result very interesting. Dominik: I am not sure we should talk about this here, as 1) I am not sure whether it is a good idea, and 2) we need a complete paragraph to make it understandable, which might be too long for the conclusion section... Following MCS's suggestion, I would rather add some more details on manifold learning.
Additionally, the investigation of more complex mechanical systems involving several members modeled with different structural elements and the extension of the current ideas to the dynamic context will be addressed in future works.

%\MCS[The sentence on manifold learning is quite vague to me.
%Is it an advanced technique?
%There are well-established mathematical methods with theoretical foundation.
%What is the embedding about?]

%% file: 7back.tex
\section*{Acknowledgments}

\noindent C. G. Gebhardt and R. Rolfes gratefully acknowledge the financial support of the Lower Saxony Ministry of Science and Culture (research project \textit{ventus efficiens}, FKZ ZN3024) and the German Research Foundation (research project ENERGIZE, GE 2773/3-1 -- RO 706/20-1) that enabled this work. D. Schillinger acknowledges support from the German Research Foundation through the DFG Emmy Noether Award SCH 1249/2-1, and from the European Research Council via the ERC Starting Grant ``ImageToSim'' (Action No.\ 759001).

%\noindent We also thank the reviewers for their valuable comments.